\documentclass[reqno,12pt]{amsart} 
\usepackage{vmargin}
\setpapersize{A4}

\usepackage{amsmath} 

\usepackage{amsfonts}   

\usepackage{amssymb}      

\usepackage{eufrak}      

\usepackage{amscd}      

\usepackage{amsthm}      
   
\usepackage{epsfig}      

\usepackage{amstext}      

\usepackage[all,line,dvips]{xy} 
\CompileMatrices 

\newcommand{\Aut}{\operatorname{Aut}}

\newcommand{\Span}{\operatorname{Span}}

\newcommand{\Per}{\operatorname{Per}}

 \newcommand{\supp}{\operatorname{supp}}

\newcommand{\Is}{\operatorname{Is}}
\newcommand{\alg}{\operatorname{alg}}
 \newcommand{\plus}{\operatorname{+}}
 \newcommand{\Prim}{\operatorname{Prim}}

   \theoremstyle{plain}
   \newtheorem{thm}{Theorem}[section]
   
   \newtheorem{lemma}[thm]{Lemma}  
   \newtheorem{cor}[thm]{Corollary}
   \theoremstyle{definition}
   
   \newtheorem{defn}[thm]{Definition}
   
   \theoremstyle{remark}
   
   \newtheorem{remark}[thm]{Remark}

\usepackage{graphicx}

   \numberwithin{equation}{section}

        \date{\today}

\title[$C^*$-algebras and groupoids]{Semi \'etale groupoids and applications}
  
\author{Klaus Thomsen}

\date{\today}

\email{matkt@imf.au.dk}
\address{Institut for matematiske fag, Ny Munkegade, 8000 Aarhus C, Denmark}

\begin{document}

\maketitle

\begin{abstract} We associate a $C^*$-algebra to a locally compact
  Hausdorff groupoid with the property that the range map is locally
  injective. The construction generalizes J. Renault's reduced
  groupoid $C^*$-algebra of an \'etale groupoid and has the advantage
  that it works for the groupoid arising from a locally injective
  dynamical system by the method introduced in increasing generality
  by Renault, Deaconu and Anantharaman-Delaroche. We study the
  $C^*$-algebras of such groupoids and give necessary and sufficient
  conditions for simplicity, and show that many of them contain a
  Cartan subalgebra as defined by Renault. In particular,
  this holds when the dynamical system is a shift space, in which
  case the $C^*$-algebra coincides with the one introduced by Matsumoto and Carlsen.  
\end{abstract}

\section{Introduction}

The main purpose of this paper is to develop new tools for the
investigation of $C^*$-algebras which have been constructed from
shift spaces in a series of papers by K. Matsumoto and
T. Carlsen, cf. \cite{Ma1}-\cite{Ma5}, \cite{C}, \cite{CM}.  The main results about the structure of these
algebras which we obtain here give necessary and sufficient
conditions for the algebras to be simple, and show that they all
contain a Cartan subalgebra in the sense introduced by
J. Renault in \cite{Re3}. Previous results on simplicity of the
$C^*$-algebras defined from subshifts are all due to Matsumoto and
give only sufficient
conditions under various additional assumptions on the
subshift. As a step on the way we show that each of these algebras is
a crossed
product in the spirit of Paschke, \cite{P}, arising from a full
corner endomorphism of an AF-algebra.

The methods we employ are useful beyond the study of
$C^*$-algebras of subshifts because they
extend the applicability of locally compact groupoids to the
construction and study of $C^*$-algebras. The use of groupoids in
relation to $C^*$-algebras was initiated by the pioneering work of J. Renault in \cite{Re1}. After a
relatively slow beginning during the eighties the last two
decades has witnessed an increasing recognition of the importance of
groupoids as a tool to encode various mathematical structures in a
$C^*$-algebra. Of particular importance in this respect are the
so-called \'etale groupoids which have been used in many different contexts, for example in
connection with graph algebras and dynamical systems. In an \'etale
groupoid the range and source maps are local homeomorphisms, and in
particular open as they must be if there is a Haar system
in the sense of Renault, cf. \cite{Re1}. But in the locally compact
groupoid which is naturally associated to a dynamical system by the
construction of Renault, Deaconu and Anantharaman-Delaroche, cf.
\cite{Re1}, \cite{De} and \cite{A}, the range and source maps are only
open if the map of the dynamical system is also open, and this is a
serious limitation which for example prevents the method from being
used on subshifts which are not of finite type. For this reason we
propose here a construction of a $C^*$-algebra from a more general class of locally compact
Hausdorff groupoids which differ from
the \'etale groupoids in that the range and source maps are locally
injective, but not necessarily open. This class of groupoids is not
new; it coincides with the locally
compact Hausdorff groupoids which were called $r$-discrete by Renault
in \cite{Re1} and they are equipped with a (continuous) Haar system if
and only if they are \'etale. In many influential places in the litterature on the
$C^*$-algebras of groupoids, such \cite{A} or \cite{Pa} for example, an r-discrete
groupoid is assumed to have a continuous Haar system and hence
to be \'etale in the terminology which is now generally
accepted. In order to avoid any misunderstanding we therefore propose
the name semi \'etale for the class of locally compact groupoids where
the range map is locally injective, but not necessarily open.

The algebra we associate to a locally
compact semi \'etale groupoid is the
reduced groupoid $C^*$-algebra of Renault when the groupoid is \'etale
and the construction is a generalization of
his. To some extend the only price one has to pay when dealing with groupoids
which are not \'etale, and only semi \'etale, is that the continuous
and compactly supported functions no longer are invariant under the
convolution product and hence do not constitute a $*$-algebra with
respect to that product. Nonetheless they still generate a $C^*$-algebra and
we obtain results on its structure which go
beyond those known in the \'etale case, as far as
necessary and sufficient conditions for simplicity and the presence of
a Cartan subalgebra is concerned. 

In the second part of the paper we make a first investigation of the $C^*$-algebras which arise from
the construction of Renault, Deaconu and Anantharaman-Delaroche when
the map of
the dynamical system is locally injective but not necessarily
open. An interesting class of such dynamical systems are the one-sided
subshifts since the shift map is locally injective but only open when
the shift space is of
finite type. We show that the (reduced) $C^*$-algebra of the
semi \'etale groupoid constructed from a one-sided subshift is a copy of the
Matsumoto-algebra of Carlsen, cf. \cite{C}, and the results concerning
its
structure are obtained by specializing results on the groupoid
$C^*$-algebra arising from a
general locally injective map. 

\bigskip

\emph{Acknowledgement:}  I want to thank Toke Meier Carlsen for
valuable information concerning the $C^*$-algebras of subshifts, and
the referee for his remarks on the first versions of this paper. It was
him (or her) who pointed out that openness
of the unit space implies local
injectivity of the range map.

\section{The $C^*$-algebra of a semi \'etale groupoid}

\subsection{Definitions and fundamental tools} \label{Sec1}

Let $G$ be a locally compact groupoid, cf. \cite{Re1}. As
in \cite{Re1} we denote the unit space of $G$ by $G^0$ and use the
letters $r$ and $s$ for the range and source maps, respectively. 
We will say that $G$ is \emph{semi \'etale} when $r : G \to G^0$ is
locally injective, i.e. when the topology of $G$
has a base consisting of open sets $U$ such that $r : U \to G^0$ and
$s : U \to G^0$ are
injective. An open subset $U \subseteq G$ with this property will be
called \emph{a bisection.}

\begin{lemma}\label{rdiscrete} $G^0$ is open in $G$ if and only if $r$
  is locally injective.
\begin{proof} Assume first that $r$ is locally injective. Let $x \in G^0$ and fix a bisection $U$ containing
  $x$. If every open neighborhood of $x$ contained
  an element from $G \backslash G^0$ the continuity of the groupoid
  operations would imply the existence of an element $\gamma \in
  U\backslash G^0$
  with $r(\gamma) \in U$. This violates the injectivity of $r$ on $U$.

Conversely, assume that $G^0$ is open in $G$. Let $\gamma \in G$. By continuity of the
groupoid operations there is an open neighborhood $U$ of $\gamma$
in $G$ such that $\mu^{-1}\nu \in G^0$ for all $\mu,\nu \in U$ with
$r(\mu) = r(\nu)$. Then $r$ is injective on $U$. 
\end{proof}
\end{lemma}

Assume now that $G$ is semi \'etale.

\begin{lemma}\label{semidiscrete} (Lemma 2.7 (i) in Chapter I of
  \cite{Re1}.)  Let $x \in G^0$. Then $r^{-1}(x)$ and $s^{-1}(x)$ are
discrete sets in the topology inherited from $G$.
\begin{proof} Let $U$ be a bisection
  containing $x$. Since $U \cap r^{-1}(x) = \{x\}$ we see that $x$
  is isolated in $r^{-1}(x)$. A similar argument shows that $y$ is
  isolated in $s^{-1}(y)$ for all $y \in G^0$. Let $\gamma \in r^{-1}(x)$ and define $\Phi : r^{-1}(x)
  \to r^{-1}\left(s(\gamma)\right)$ such that $\Phi(\eta) =
  \gamma^{-1}\eta$. Then $\Phi$ is a homeomorphism with inverse $\eta
  \mapsto \gamma \eta$. Since $\Phi(\gamma) = s(\gamma)$ and $s(\gamma)$
  is isolated in $r^{-1}\left(s(\gamma)\right)$ it follows that
  $\gamma$ is isolated in $r^{-1}(x)$. This proves that $r^{-1}(x)$ is discrete. The argument concerning $s^{-1}(x)$ is identical.
\end{proof}
\end{lemma}

It follows from Lemma \ref{semidiscrete} that $r^{-1}(x)$ and
$s^{-1}(x)$  both have a finite intersection
with any compact subset of $G$. Therefore, when $f,g : G \to \mathbb
C$ are compactly supported functions, we can define $f \star g: G \to
\mathbb C$ by the usual formula
\begin{equation}\label{stareq}
f \star g(\gamma) = \sum_{\gamma_1 \gamma_2 = \gamma}
f(\gamma_1)g(\gamma_2) .
\end{equation}
Then $f \star g$ is again compactly supported, and $f\star g$ is bounded when $f$ and
$g$ both are. It follows that the set $B_c(G)$ of bounded compactly supported
functions on $G$ is a $*$-algebra with the product $\star$ and the
involution $f \mapsto f^*$ defined such that
\begin{equation}\label{stareq2}
f^*(\gamma) = \overline{f\left(\gamma^{-1}\right)} .
\end{equation}
To obtain a $C^*$-norm we use the usual representations: For each $x
\in G^0$ we define a $*$-representation $\pi_x$ of $B_c(G)$ on
$l^2\left(s^{-1}(x)\right)$ such that
$$
\left(\pi_x(f)\psi\right)(\gamma) =  \sum_{ \gamma_1 \gamma_2 = \gamma}
\ f\left(\gamma_1\right)\psi\left(\gamma_2\right) .
$$ 
We define the $C^*$-algebra $B^*_r(G)$ to be the completion of
$B_c(G)$ in the norm
$$
\left\|f\right\| = \sup_{x \in G^0} \left\|\pi_x(f)\right\| .
$$ 
Let $C_c(G)$ be the subspace of $B_c(G)$ consisting of the functions
on $G$ which are compactly supported and continuous. We let $C^*_r(G)$ be the $C^*$-subalgebra of $B^*_r(G)$ generated
by $C_c(G) \subseteq B_c(G)$, i.e.
$$
C^*_r(G) = \overline{\alg^* G}
$$
when $\alg^* G$ denotes the $*$-subalgebra of $B^*_r(G)$ generated
by $C_c(G)$. Note that $C^*_r(G)$ is separable when
$G$ is second countable while $B^*_r(G)$ essentially never is.

\begin{lemma}\label{estimates} (Proposition 4.1 in Chapter II of \cite{Re1}.)
Let $f \in B_c(G)$. Then 
\begin{equation}\label{on1}
\sup_{\gamma \in G} \left|f(\gamma)\right| \leq \|f \|
\end{equation}
and 
\begin{equation}\label{on2}
\sum_{\gamma \in s^{-1}(x)} |f(\gamma)|^2 \leq \|f\|^2 
\end{equation}
for all $x \in G^0$.
\begin{proof} Let $\gamma \in G$ and set $x = s(\gamma)$. Let $\delta_x,
  \delta_{\gamma} \in l^2\left(s^{-1}(x)\right)$ denote the
  characteristic functions of $\left\{x\right\}$ and
  $\left\{\gamma\right\}$, respectively. Then
  $\left(\pi_x(f)\delta_x\right)(\gamma) = \sum_{\gamma_1 \gamma_2 = \gamma}
f\left(\gamma_1\right)\delta_x\left(\gamma_2\right) =
f\left(\gamma\right)$ which shows that
\begin{equation}\label{eq21}
\left<\pi_x(f)\delta_x,\delta_{\gamma}\right> = f\left(\gamma\right),
\end{equation}
and
\begin{equation}\label{eq123}
\left\|\pi_x(f)\delta_x\right\|^2_{l^2\left(s^{-1}(x)\right)} =
\sum_{\gamma \in s^{-1}(x)} \left|f(\gamma)\right|^2 .
\end{equation}
(\ref{on1}) follows from (\ref{eq21}) and (\ref{on2}) follows from (\ref{eq123}).
\end{proof}
\end{lemma}

\begin{lemma}\label{renaultineq} Let $f \in B_c(G)$ be supported in a
  bisection. Then $\left\|f\right\| = \sup_{\gamma \in
    G} \left|f(\gamma)\right|$.
\begin{proof} Let $U$ be a bisection containing $\supp f$. Define $\tilde{f} : G \to \mathbb C$ such that
  $\tilde{f} (\gamma) = 0$ when $r(\gamma) \notin r(U)$ and
  $\tilde{f}(\gamma) = \overline{f( \mu)}$ where $\mu \in U$ is the unique element
  with $r(\mu) = r(\gamma)$ when $r(\gamma) \in r(U)$. Let $x \in G^0$
  and define $V : l^2(s^{-1}(x)) \to l^2(s^{-1}(x))$ such that
$V\varphi(\gamma) = 0$ when $r(\gamma) \notin
    r(U)$ and $V\varphi(\gamma) = \varphi(\mu^{-1}\gamma)$ when $r(\gamma) \in r(U)$, where $\mu \in U$ is the element
  with $r(\mu) = r(\gamma)$. Then $\left\|V\right\| \leq 1$. Let $\varphi, \psi \in l^2\left(s^{-1}(x)\right)$. Then
\begin{equation*}
\begin{split} 
& \left|\left< \pi_x(f) \varphi, \psi \right> \right| = \left|  \sum_{\gamma \in s^{-1}(x)} \sum_{\gamma_1 \gamma_2 = \gamma} f(\gamma_1)\varphi(\gamma_2)\overline{\psi(\gamma)}\right| \\
& = \left|\sum_{\gamma \in s^{-1}(x)} V\varphi(\gamma)
\overline{\tilde{f}(\gamma)\psi(\gamma)}  \right|
 \leq \left\|V\varphi\right\| \left\|\psi\right\| \sup_{\gamma} \left|\tilde{f}(\gamma)\right| \leq  \|\varphi\| \|\psi\| \sup_{\gamma \in G} \left|f(\gamma)\right|.
\end{split}
\end{equation*}
It follows that $\left\|f\right\| \leq \sup_{\gamma \in
    G} \left|f(\gamma)\right|$. Equality holds by (\ref{on1}).
\end{proof}
\end{lemma}

Let $B_0(G)$ denote the space of bounded functions on $G$ which
vanishes at infinity. We consider $B_0(G)$ as a Banach space in the
supremum norm $\left\| \cdot \right\|_{\infty}$. It follows from (\ref{on1}) that the inclusion
$B_c(G) \subseteq B_0(G)$ extends to a linear map $j : B^*_r(G) \to
B_0(G)$ such that
\begin{equation}\label{kasper}
\left\|j(b)\right\|_{\infty} \leq \left\|b\right\|_{B^*_r(G)} .
\end{equation}

\begin{lemma}\label{onj} (Proposition 4.2 (iii) in Chapter 3 of \cite{Re1}.) Let $a,b \in B_r^*(G)$. Then $j(b)|_{s^{-1}(x)}
  \in l^2\left(s^{-1}(x)\right),
  \ j(a)|_{r^{-1}(x)} \in l^2\left(r^{-1}(x)\right)$ for all $x \in G^0$, and
\begin{equation}\label{on17}
j(a^*)(\gamma) = \overline{j(a)(\gamma^{-1})} ,
\end{equation}
and
\begin{equation}\label{on4}
j(ab)(\gamma) = \sum_{\gamma_1\gamma_2 = \gamma}
j(a)\left(\gamma_1\right)j(b)\left(\gamma_2\right)
\end{equation}
for all $\gamma \in G$.
\begin{proof} Choose sequences $\{f_n\}, \left\{g_n\right\} \subseteq
  B_c(G)$ such that $a = \lim_{n \to \infty} f_n$ and $b = \lim_{n \to
    \infty} g_n$ in $B^*_r(G)$. It follows from (\ref{on1}) that
  
$$
j(a^*)(\gamma) = \lim_{n \to \infty} f_n^*(\gamma) = \lim_{n \to
    \infty} \overline{f_n(\gamma^{-1})} =
  \overline{j(a)\left(\gamma^{-1}\right)}
$$ 
which gives (\ref{on17}). It follows from (\ref{on2}) that $j(b)|_{s^{-1}(x)}$
is the limit in $l^2\left(s^{-1}(x)\right)$ of the sequence
$\left\{g_n|_{s^{-1}(x)} \right\}$. Inserting $f^*$ for $f$ in
(\ref{on2}) we obtain the inequality    
\begin{equation}\label{on3}
\sum_{\gamma \in r^{-1}(x)} \left|f(\gamma)\right|^2 \leq \|f\|^2
\end{equation}
when $x \in G^0$ and $f \in B_c(G)$. Then (\ref{on3}) implies that $j(a)|_{r^{-1}(x)}$
is the limit in $l^2\left(r^{-1}(x)\right)$ of the sequence
$\left\{f_n|_{r^{-1}(x)} \right\}$. In particular, $j(b)|_{s^{-1}(x)}$
and $ j(a)|_{r^{-1}(x)}$ are both square-summable functions for all $x
\in G^0$ and hence the righthand side of (\ref{on4}) makes sense for
each $\gamma \in G$. Let $\gamma \in G$ and set $x = r(\gamma), \ y =
s(\gamma)$. We have then the estimate
\begin{equation*}
\begin{split}
&  \left| \sum_{\gamma_1\gamma_2 = \gamma}
  j(a)\left(\gamma_1\right)j(b)\left(\gamma_2\right)  -
  j(f_n \star g_n)(\gamma)\right| \\
& =\left| \sum_{\gamma_1\gamma_2 = \gamma}
  j(a)\left(\gamma_1\right)j(b)\left(\gamma_2\right)  -
  \sum_{\gamma_1\gamma_2 = \gamma}
  f_n\left(\gamma_1\right)g_n\left(\gamma_2\right) \right|  \\
& \leq \sum_{\gamma_1\gamma_2 = \gamma} 
\left|j(a)\left(\gamma_1\right) - f_n\left(\gamma_1\right)\right|\left|j(b)\left(\gamma_2\right)\right| + \sum_{\gamma_1\gamma_2 = \gamma}
\left|f_n\left(\gamma_1\right)\right|\left|g_n\left(\gamma_2\right) -
  j(b)\left(\gamma_2\right)\right| \\
& \leq \left\|j(a)
  -f_n\right\|_{l^2\left(r^{-1}\left(x\right)\right)}
\left\|j(b)\right\|_{l^2\left(s^{-1}\left(y\right)\right)} +
\left\|f_n\right\|_{l^2\left(r^{-1}\left(x\right)\right)}
\left\|g_n -
  j(b)\right\|_{l^2\left(s^{-1}\left(y\right)\right)} .
\end{split}
\end{equation*}
The equality (\ref{on4}) follows then by letting $n$ tend to infinity.
\end{proof}
\end{lemma}

\begin{cor}\label{jinj} (Proposition 4.2 (i) in Chapter 3 of \cite{Re1}.) $j : B^*_r(G) \to B_0(G)$ is injective.
\begin{proof} If $j(a) = 0$ it follows from (\ref{on4}) that $j(ab) = 0$ for all $b \in B^*_r(G)$. Now note that
it follows from (\ref{on1}) that the equality (\ref{eq21}) extends by
continuity to the equality
$$
\left<\pi_x(d)\delta_x, \delta_{\gamma} \right> =
j(d)\left(\gamma\right),
$$
valid for all $d \in B^*_r(G)$, all $x \in G^0$ and all $\gamma\in s^{-1}(x)$. Since $j(ab) =
0$ for all $b \in B^*_r(G)$ this implies that
$$
\left<\pi_x(a)\pi_x(b)\delta_x, \delta_{\gamma} \right> = \left<\pi_x(ab)\delta_x, \delta_{\gamma} \right> = 0
$$
for all $x \in G^0, b \in B^*_r(G)$ and all $\gamma\in
s^{-1}(x)$. Since $\delta_x $ is cyclic for $\pi_x$ this implies that
$\pi_x(a) = 0$ for all $x$, i.e. $a = 0$.
\end{proof}
\end{cor}

Since $G^0$ is closed in $G$ we have an embedding
$B_c(G^0) \subseteq B_c(G)$. Let $B_0(G^0)$ be the $C^*$-algebra of
bounded functions on $G^0$ which vanish at infinity. Note that $\sup_{x
  \in G^0}\left\|\pi_x(f)\right\| = \sup_{y \in G^0}
\left|f(y)\right|$ when $f \in B_c(G^0)$ by Lemma
\ref{renaultineq}. It follows that the embedding $B_c(G^0) \subseteq B_c(G)$ extends
by continuity to an isometric $*$-homomorphism $B_0(G^0)\to
B^*_r(G)$. In the following we will consider $B_0(G^0)$ as a
$C^*$-subalgebra of $B^*_r(G)$ via this embedding. It follows from
(\ref{on1}) and (\ref{eq21}) that there is a conditional expectation
$$
P_G : B^*_r(G) \to B_0(G^0)
$$
defined such that $P_G(a)(x) = \left< {\pi}_x(a)\delta_x,\delta_x
\right>$. Then
\begin{equation}\label{P=j}
P_G(a)(x) = j(a)(x)
\end{equation}
for all $a \in B^*_r(G)$, $x \in G^0$.

\begin{lemma}\label{techlemma} Let $E \subseteq
  G$ be a subset which is both closed and open in $G$. Let $f_1,f_2,
  \dots, f_n \in C_c(G)$, and let $V_{\alpha}, \alpha \in I$, be a
  collection of open sets in $G$ such that $E \subseteq
  \bigcup_{\alpha \in I} V_{\alpha}$. 

It follows that there are functions $h^j_1,h^j_2, \dots, h^j_n \in
C_c(G), \ j =1,2, \dots,m$, such that
\begin{enumerate}
\item[a)] $\sum_{j=1}^m h^j_1 \star h^j_2 \star \dots \star h^j_n
  (\gamma) = \begin{cases} f_1 \star f_2 \star \dots \star f_n(\gamma)
    , & \ \gamma \in E \\ 0, & \ \gamma \notin E,\end{cases}$ 
\end{enumerate}

and

\begin{enumerate}
\item[b)] for each $j \in \{1,2, \dots, m\}$ there is an $\alpha_j \in
  I$ such that 
$$\supp h^j_1 \star h^j_2 \star \dots \star h^j_n \
  \subseteq \ V_{\alpha_j}.
$$
\end{enumerate}
\begin{proof} We say that a function $k : G^n \to \mathbb C$ is of
  \emph{product type} when there are functions $k_1,k_2, \dots, k_n
  \in C_c(G)$ such that
$$
k(\gamma_1, \gamma_2, \dots, \gamma_n) = k_1(\gamma_1)k_2(\gamma_2)
\dots k_n(\gamma_n)
$$
for all $(\gamma_1, \gamma_2, \dots, \gamma_n) \in G^n$. Set 
$$
G^{(n)} = \left\{ (\gamma_1, \gamma_2, \dots, \gamma_n) \in G^n : \
  s\left(\gamma_i\right) = r\left(\gamma_{i+1}\right), \ i = 1,2,
  \dots, n-1\right\} .
$$
For each $\alpha \in I$, set
$$
A_{\alpha} = \left\{ \left( \gamma_1, \gamma_2, \dots, \gamma_n\right)
  \in G^{(n)}: \ \gamma_1\gamma_2 \dots \gamma_n \in V_{\alpha} \cap
  E \right\} 
$$
which is an open subset of $G^{(n)}$. Let
$$
A = \left\{ \left( \gamma_1, \gamma_2, \dots, \gamma_n\right)
  \in G^{(n)}: \ \gamma_1\gamma_2 \dots \gamma_n \in 
  E \right\} 
$$
and note that $A$ is both open and closed in $G^{(n)}$.
Let $\Omega_{\alpha} \subseteq
G^n$ be an open subset such that $\Omega_{\alpha} \cap G^{(n)} =
A_{\alpha}$. Since $\left(\supp f_1 \times \supp f_2 \times \dots \times
\supp f_n\right) \cap A \cap G^{(n)}$ is a compact subset of $G^n$
contained in $\bigcup_{\alpha \in I} \Omega_{\alpha}$ there is a
cover $\Omega'_{\beta}, \beta \in I'$, of 
$$
\left(\supp f_1 \times \supp f_2 \times \dots \times
\supp f_n\right) \cap A \cap G^{(n)}
$$ 
in $G^n$ such that each $\Omega'_{\beta}$
is an open rectangle, i.e. of the form
$$
\Omega'_{\beta} = U_1 \times U_2 \times \dots \times U_n,
$$
where each $U_i$ is an open subset of $G$, and such that the closure,
$\overline{\Omega'_{\beta}}$, of each $\Omega'_{\beta}$ is contained
in $\Omega_{\alpha}$ for some $\alpha$. By compactness there is a
finite set $\left\{\beta_1,\beta_2, \dots , \beta_{m'}\right\} \subseteq I'$ such that
$$
\left(\supp f_1 \times \supp f_2 \times \dots \times
\supp f_n\right) \cap A \cap G^{(n)} \subseteq \bigcup_{j = 1}^{m'}
\Omega'_{\beta_j} .
$$   
For each $j \in \{1,2, \dots, m'\}$ there is a positive function $g_j
\in C_c\left(G^n\right)$ of product type such that $g_j(\xi) = 1, 
\xi \in \overline{\Omega'_{\beta_j}}$, and $\supp g_j \subseteq
\Omega_{\alpha_j}$ for some $\alpha_j \in I$ with
$\overline{\Omega'_{\beta_j}} \subseteq \Omega_{\alpha_j}$. Define
$h_j, j = 1,2, \dots, m'$, such that $h_1 = g_1$ and $h_{i+1} = (1 -
g_1)(1-g_2) \dots (1 -g_i)g_{i+1}, 1\leq i \leq m'-1$. Then $h_1(\xi)
+ h_2(\xi) + \dots + h_{m'}(\xi) = 1$ when $\xi \in \left(\supp f_1 \times \supp f_2 \times \dots \times
\supp f_n\right) \cap A \cap G^{(n)}$. Furthermore, each $h_j$ is the
sum of functions of product type, each of which has its support
contained in some $\Omega_{\alpha}$. Let $h'_j, j =1,2, \dots, m$, be
an enumeration of these functions such that $\sum_{j=1}^m h'_j =
\sum_{j=1}^{m'} h_j$. Since $h'_j$ is of product type there are functions
$k^j_1,k^j_2, \dots, k^j_n \in C_c(G)$ such that
$$
h'_j(\gamma_1, \gamma_2, \dots, \gamma_n) =  k^j_1(\gamma_1)k^j_2(\gamma_2)
\dots k^j_n(\gamma_n)
$$
for all $(\gamma_1, \gamma_2, \dots, \gamma_n) \in G^n$. Set $h^j_i =
k^j_if_i$. Then $h^j_1,h^j_2, \dots, h^j_n \in
C_c(G), \ j =1,2, \dots,m$, satisfy a)
and b) by construction. 
\end{proof}
\end{lemma}

\begin{lemma}\label{condexp} 
\begin{enumerate}
\item[i)] $P_G$ is positive, i.e. $a \geq 0$ in $B^*_r(G)$ $\Rightarrow$
  $P_G(a) \geq 0$ in $B_0\left(G^0\right)$.
\item[ii)] $P_G(b) = b$ when $b \in B_0\left(G^0\right)$. 
\item[iii)] $\|P_G\| = 1$.
\item[iv)] $P_G$ is faithful, i.e. $a \neq 0 \ \Rightarrow \ P_G(a^*a)
  \neq 0$.
\item[v)] $P_G\left(C^*_r(G)\right) = C^*_r(G) \cap
  B_0(G^0) = \overline{\alg^* G \cap B_c(G^0)}$.
\end{enumerate}
\begin{proof} i), ii) and iii) hold by construction.

iv): Let $a
\in B^*_r(G)$ and assume that $P_G(a^*a) = 0$. It follows then from
(\ref{P=j}) and (\ref{on4}) that
$$
\sum_{\mu \in s^{-1}(x)} \left|j(a)\left(\mu\right)\right|^2  =\sum_{\gamma_1 \gamma_2 =x}
\overline{j(a)\left(\gamma_1^{-1}\right)}j(a)\left(\gamma_2\right) =
j(a^*a)(x) =0
$$
for all $x \in G^0$. This shows that $j(a) = 0$ and it follows then from
Corollary \ref{jinj} that $a = 0$.

v) : The inclusions $\overline{\alg^* G \cap B_c(G^0)} \subseteq
C^*_r(G) \cap B_0(G^0) \subseteq P_G\left(C^*_r(G)\right)$ are obvious
so it suffices to show that $P_G\left(\alg^* G\right) \subseteq
\alg^* G \cap B_c(G^0)$. Since $P_G(a) = j(a)|_{G^0} = a|_{G^0}$ when
$a \in \alg^* G$, this follows from Lemma \ref{techlemma}, applied
with $E = G^0$.

\end{proof}
\end{lemma}

It follows from Lemma \ref{condexp} and a result of Tomiyama that 
\begin{equation}\label{tomiyama}
P_G(d_1ad_2) = d_1P_G(a)d_2
\end{equation}
for all $a \in C^*_r(G), \ d_1,d_2 \in C^*_r(G) \cap
  B_0(G^0)$; a fact which can also easily be established directly.

\begin{lemma}\label{type1} Assume that $n \in \alg^*G$ is supported in
  a bisection. It follows that
\begin{equation}\label{nPcomm}
n^*P_G(a)n = P_G(n^*an)
\end{equation}
for all $a \in C^*_r(G)$.
\begin{proof} Let $U$ be a bisection containing $\supp n$. It follows
  from Lemma \ref{onj} that
\begin{equation*}
\begin{split}
&j(n^*P_G(a)n)(\gamma) = \sum_{\gamma_1\gamma_2\gamma_3 = \gamma}
\overline{n(\gamma_1^{-1})}j(P_G(a))(\gamma_2)n(\gamma_3) \\
& = \begin{cases} 0, & \ \text{when} \ \gamma \notin s(U) \\
  \overline{n(\mu)} j(a)(r(\mu)) n(\mu) \ & \text{where $\mu\in U \cap
    s^{-1}(\gamma)$, when $\gamma \in s(U)$ .}
\end{cases}  
\end{split}
\end{equation*}
This is the same expression we find for $j(P_G(n^*an))(\gamma)$ and
hence (\ref{nPcomm}) follows from Corollary \ref{jinj}.
\end{proof}
\end{lemma}

\begin{lemma}\label{openincl} Let $H \subseteq G$ be an open
  subgroupoid, i.e. $H$ is open, $H^{-1} = H$ and $\left(\gamma_1,\gamma_2\right) \in H^2 \cap
  G^{(2)} \Rightarrow \gamma_1 \gamma_2 \in H$. Then the inclusions $C_c(H) \subseteq
  C_c(G)$ and $B_c(H) \subseteq B_c(G)$ extend to $C^*$-algebra
  embeddings $C^*_r(H) \subseteq
  C^*_r(G)$ and $B^*_r(H) \subseteq B^*_r(G)$, respectively.
\begin{proof} Clearly, the inclusion $C_c(H) \subseteq
  C_c(G)$ extends to an inclusion $\alg^* H \subseteq \alg^* G$ of
  $*$-algebras so
  it remains only to show that $\left\|f\right\|_{B^*_r(G)} =
  \left\|f\right\|_{B^*_r(H)}$ when $f \in B_c(H)$. To this end,
  let $x \in G^0$ and let $f \in B_c(H)$. We define an equivalence relation $\sim$ on
  $Hs^{-1}(x)$ such that $\gamma \sim \gamma'$ if and only if $\gamma'
  = \mu \gamma$ for some $\mu \in H$. Let $\left[Hs^{-1}(x)\right]$
  denote the set of equivalence classes in $Hs^{-1}(x)$. Then
$$
l^2\left(s^{-1}(x)\right) = \oplus_{\xi \in \left[Hs^{-1}(x)\right]}
l^2(\xi) \oplus l^2\left(s^{-1}(x) \backslash Hs^{-1}(x)\right) 
$$
and $\pi_x(f)$ respects this direct sum decomposition. Since $\pi_x(f)
= 0$ on $l^2\left(s^{-1}(x) \backslash Hs^{-1}(x)\right)$ 
we find that
\begin{equation}\label{est1}
\left\|\pi_x(f)\right\| = \sup_{\xi \in \left[Hs^{-1}(x)\right]}
\left\|\pi_x(f)|_{l^2(\xi)}\right\| .
\end{equation}
Let $\xi \in \left[Hs^{-1}(x)\right]$ and fix a representative $\gamma_0 \in
\xi$. We can then define a unitary $V : l^2\left(H \cap
  s^{-1}\left(r(\gamma_0)\right)\right) \to l^2(\xi)$ such that
$V\psi(\mu) = \psi\left(\mu\gamma_0^{-1}\right)$.
It is then straightforward to verify that
$$
V\pi_{r(\gamma_0)}(f)|_{l^2\left(H\cap s^{-1}(r(\gamma_0))\right)}V^*
= \pi_x(f)|_{l^2(\xi)} ,
$$
and we conclude that 
$$
\left\|\pi_x(f)|_{l^2(\xi)} \right\| = \left\|
  \pi_{r(\gamma_0)}(f)|_{l^2\left(H\cap
      s^{-1}(r(\gamma_0))\right)}\right\| \leq
\left\|f\right\|_{B^*_r(H)} .
$$
Combined with (\ref{est1}), and using that $x \in G^0$ was arbitrary, this shows that
$\left\|f\right\|_{B^*_r(G)} \leq \left\|f\right\|_{B^*_r(H)}$. Since
the reversed inequality is trivial, this completes the proof.  
\end{proof}
\end{lemma}

\begin{lemma}\label{exelogrenault} Let $g \in C_c(G)$ and $f \in
  \alg^* G$. Then the pointwise product
$$
g \cdot f(\gamma) = g(\gamma)f(\gamma), \ \ \gamma \in G,
$$
is in $\alg^*G$.
\begin{proof} It follows from Lemma \ref{techlemma}, applied with $E =
  G$, that $f$ is a finite sum of elements from $\alg^*G$ whose
  compact supports are contained in bisections. We may therefore
  assume that $f$ has support in a bisection. Then an argument from
  Lemma 4.3 of \cite{ER} completes the proof: Define first $u_0 : r(\supp f) \to \mathbb
  C$ such that $u_0(x) = g(\mu)$, where $\mu \in \supp f$ is the
  unique element with $r(\mu) = x$, and let $u \in C_c(G^0)$ be an
  extension of $u_0$. Then $g(\gamma)f(\gamma) = u(r(\gamma))f(\gamma)
  = u \star f(\gamma)$ for all $\gamma \in G$.
\end{proof}
\end{lemma}

\begin{lemma}\label{product2} Let $a \in C^*_r(G)$ and let $h \in
 C_c(G)$ be supported in a bisection.
There is then an element
  $h \cdot a \in C_r^*(G)$ such that $j(h \cdot a)(\gamma) =
  h(\gamma)j(a)(\gamma)$ for all $\gamma \in G$.
\begin{proof} Define a function $\tilde{h} : G \to \mathbb C$ such
  that $
\tilde{h}(\gamma) = $0 when $r(\gamma) \notin
  s(\supp h)$ and $\tilde{h}(\gamma) = h(\gamma')$ where $\gamma' \in \supp h$ is
      the unique element of $s^{-1}(r(\gamma)) \cap \supp h$ when 
      $r(\gamma) \in s(\supp h)$. Then
$$
\sum_{\gamma_1 \gamma_2 = \gamma}
h(\gamma_1)f(\gamma_1)\varphi(\gamma_2) = \sum_{\gamma_1 \gamma_2 = \gamma}
f(\gamma_1)\tilde{h}(\gamma_2)\varphi(\gamma_2)
$$
when $f \in \alg^* G$ and $\varphi \in
l^2\left(s^{-1}(x)\right)$. It follows that $\left\| h \cdot f\right\|
\leq \|f\| \left\|h\right\|_{\infty} $ for all $f \in
\alg^* G$. In particular, it follows that $\{h\cdot f_n\}$ converges
in $B^*_r(G)$ when $\{f_n\} \subseteq \alg^* G$ converges to
$a$. It follows from Lemma \ref{exelogrenault} that the limit $h \cdot a =
\lim_{n \to \infty} h\cdot f_n$ exists in $C^*_r(G)$. The limit will have
the stated property since $j(h \cdot a)(\gamma) = \lim_{n \to \infty} j(h\cdot
f_n)(\gamma)$ for all $\gamma$.
\end{proof}
\end{lemma}

\subsection{Ideals}
 
Let $A$ be a
$C^*$-algebra and $D \subseteq A$ an abelian $C^*$-subalgebra. An
element $a\in A$ is a \emph{$D$-normalizer} when $a^*Da \subseteq D$
and $aDa^* \subseteq D$. The set of $D$-normalizers will be denoted by
$N(D)$.

Consider now the case where $A = C^*_r(G)$ 
and 
$$
D = D_G =C^*_r(G) \cap
  B_0(G^0).
$$

Let $N_0(D_G)$ denote the set of functions $g$ from $\alg^* G$ that
are supported in a bisection. It follows from Lemma \ref{type1} that $N_0(D_G)
\subseteq N(D_G)$. A (closed) ideal $J \subseteq D_G$ is said to be \emph{$G$-invariant} when $n^*Jn
\subseteq J$ for all $n \in N_{0}(D_G)$. Note that $I\cap D_G$ is a
 $G$-invariant ideal in $D_G$ when $I$ is a (closed and twosided) ideal in $C^*_r(G)$.

\begin{lemma}\label{invideal} Let $J \subseteq D_G$ be a $G$-invariant
  ideal. It follows that
$$
\widehat{J} = \left\{ a \in C^*_r(G): \ P_G(a^*a) \in J\right\} 
$$
is an ideal in $C^*_r(G)$ such that $J = \widehat{J} \cap D_G$.
\begin{proof} It follows easily, by using the relations $x^*y^*yx \leq \|y\|^2x^*x$
and $(x+y)^*(x+y) \leq 2x^*x + 2y^*y$,
that $\widehat{J}$ is a left ideal in $C^*_r(G)$. It follows from Lemma
\ref{type1} that $a \in
\widehat{J} \ \Rightarrow \ an \in \widehat{J}$ when $n \in
N_{0}(D)$ because $J$ is $G$-invariant. It follows from Lemma
\ref{exelogrenault} that the elements of $N_{0}(D)$ span a dense
subspace in $C^*_r(G)$. We conclude therefore that $\widehat{J}$ is also a right-ideal. This proves the
lemma because the identity $J = \widehat{J} \cap D_G$ is obvious.   
\end{proof}
\end{lemma}

Note that it follows from Lemma \ref{invideal} that the lattice of
$G$-invariant ideals in $D_G$ has a copy inside the lattice of ideals in
$C^*_r(G)$.

An ideal in a $C^*$-algebra is said to be \emph{non-trivial} when it
is neither $\{0\}$ nor the whole algebra. With this terminology we
have

\begin{cor}\label{idealcor1} Assume that $D_G$ contains a non-trivial
  ideal which is $G$-invariant. It follows that $C^*_r(G)$ contains a
  non-trivial ideal.
\end{cor}

For $x \in G^0$ we let $G_x = \left\{\gamma \in G: \ r(\gamma) =
  s(\gamma) = x\right\}$ denote the \emph{isotropy group at $x$.}

\begin{lemma}\label{idealsny1} Let $I \subseteq C^*_r(G)$ be an ideal
  such that $I \cap D_G = \{0\}$. It follows that
$$
j(a)(x) = 0
$$
for all $a \in I$ and all $x \in G^0$ with $G_x = \{x\}$.
\begin{proof} Let $h \in \alg^* G$ and let $x \in G^0$ be a point with
  trivial isotropy (i.e. $G_x = \{x\})$. We assume that $h(x) \neq
  0$. Consider a point $\gamma \in G$. If $r(\gamma) = x$ and $\gamma \neq x$, we know that
  $s(\gamma) \neq x$. There is therefore an open neighborhood
  $U_{\gamma}$ of $\gamma$ such that 
  $r\left(\overline{U_{\gamma}}\right) \cap
  s\left(\overline{U_{\gamma}}\right) = \emptyset$. If $r(\gamma) \neq
  x$ there is an open neighborhood $U_{\gamma}$ of $\gamma$ such that
  $x \notin r\left(\overline{U_{\gamma}}\right)$. Finally, if $\gamma = x$
  there is an open neighborhood $U_{\gamma}$ of $\gamma$ such that
  $U_{\gamma} \subseteq G^0$. It follows from Lemma \ref{techlemma},
  applied with $E =G$,
  that there are elements $h_i \in \alg^* G$ and distinct elements $\gamma_i \in G$
  such that $\supp h_i \subseteq U_{\gamma_i}, \ i = 1,2, \dots, N$,
  and $h = \sum_{i=1}^N h_i$. By construction $x$ is only
  element of one member from $U_{\gamma_i}, i = 1,2, \dots, N$. For convenience we assume
  that $x \in U_{\gamma_1}$. Then $\gamma_1 = x$ and
  $U_{\gamma_1}\subseteq G^0$. For each $ j \geq 2$, $ x \notin
  r\left(\overline{U_{\gamma_j}}\right)$ or $x\notin
    s\left(\overline{U_{\gamma_j}}\right)$. There is therefore a function $f \in C_c(G^0)$ such
  that $0 \leq f \leq 1, \ f(x) = 1$ and $f \star h_j =0$ or $h_j
  \star f = 0, \ j \geq 2$.
 It
follows that $
f \star h_i \star f = 0$ when $i \neq 1$. Hence $f \star h \star f= f \star h_1
\star f \in D_G$. Let $q :
C^*_r(G) \to C^*_r(G)/I$ be the quotient map. Since $q$ is injective on
$D_G$ we find that
\begin{equation*}
\begin{split}
& \left\|q(h)\right\| \geq \left\|q(f \star h \star f)\right\|
= \left\|q\left(f \star h_1 \star f\right)\right\| \\
&=
\left\|f \star h_1 \star f\right\| = sup_{y\in G^0} \left|f
  \star h_1 \star f(y)\right| \geq |h(x)| .
\end{split}
\end{equation*} 
Let $a\in C^*_r(G)$. There is a sequence $\{h_k\} \subseteq \alg^*
G$ such that $a = \lim_{k \to \infty} h_k$ in $C^*_r(G)$. It follows
that
$$
\left\|q(a)\right\| =\lim_{k \to \infty} \left\|q(h_k)\right\| \geq
\lim_{k \to \infty} \left|h_k(x)\right| = \left|j(a)(x)\right|.
$$
This proves the lemma.
\end{proof}
\end{lemma}

\begin{lemma}\label{ideal13} Assume that $G_x = \{x\}$ for some $x \in
  G^0$. Let $I$ be a non-trivial ideal in $C^*_r(G)$. It
  follows that either $I\cap D_G$ or $\overline{P_G(I)}$ is a non-trivial
  $G$-invariant ideal in $D_G$.
\begin{proof} Unless the intersection $I\cap D_G$ is zero it will
  constitute an ideal in $D_G$ which must be non-trivial because $D_G$ contains an
  approximate unit for $C^*_r(G)$. Since $I\cap D_G$ is $G$-invariant it
  suffices to show that $\overline{P_G(I)}$ is a non-trivial
  $G$-invariant ideal in $D_G$ when $I \cap D_G = \{0\}$. First
  observe that it is an ideal because of (\ref{tomiyama}). Since $P_G$ is faithful by
  iv) of Lemma \ref{condexp} we have that $P_G(I) \neq 0$ since $I \neq
  0$. By assumption there is a point $x \in G^0$ with trivial isotropy
  and it follows then from Lemma \ref{idealsny1} and (\ref{P=j}) that
  $g(x) = 0$ for all $g \in \overline{P_G(I)}$. In particular,
  $\overline{P_G(I)} \neq D_G$. Thus $\overline{P_G(I)}$ is a non-trivial
  ideal in $D_G$ when $I \cap D_G$ fails to be. It is $G$-invariant by
  Lemma \ref{type1}.
\end{proof}
\end{lemma}

\begin{thm}\label{simplicity1}  Assume that $G_x = \{x\}$ for some $x \in
  G^0$. Then $C^*_r(G)$ is simple if and only if there
  are no non-trivial $G$-invariant ideals in $D_G$.
\begin{proof} Combine Lemma
  \ref{ideal13} and Corollary \ref{idealcor1}.
\end{proof}
\end{thm}

For the formulation of the following corollary remember that a subset
$V \subseteq G^0$ is $G$-invariant when $\gamma \in G, \ s(\gamma) \in
V \ \Rightarrow \ r(\gamma) \in V$.

\begin{cor}\label{invopne} Assume that $G$ is \'etale and that $G_x
  =\{x\}$ for some $x \in G^0$. It follows that $C^*_r(G)$ is simple
  if and only if there are no open non-trivial $G$-invariant subset of
  $G^0$.
\begin{proof} Since $G$ is \'etale, $D_G = C_0(G^0)$. Let $U \subseteq
  G^0$ be an open subset. By Theorem \ref{simplicity1} it suffices to
  show that the ideal $C_0(U)$ of $D_G$ is $G$-invariant if and only $U$
  is $G$-invariant. Assume first that $C_0(U)$ is $G$-invariant and let $\gamma \in G$ be such that
  $s(\gamma) \in U$. There is then an element $h \in N_0(D_G)$ such that
  $h(\gamma) =1$. It follows that $h^*h(s(\gamma)) =
  \left|h(\gamma)\right|^2 = hh^*(r(\gamma)) =1$. Since $s(\gamma) \in
  U$ there is an $f \in C_0(U)$ such that $f^*h^*hf \in C_0(U)$ and
  $f^*h^*hf(s(\gamma)) = 1$. Since $hf \in N_0(D_G)$ we find that
  $hf\left(f^*h^*hf\right)f^*h^* \in C_0(U)$ and hence that $hff^*h^*
  \in C_0(U)$. Since $hff^*h^*(r(\gamma)) = f^*h^*hf(s(\gamma)) = 1$
  this implies that $r(\gamma) \in U$.

Assume next that $U$ is $G$-invariant and let $f \in C_0(U), h \in
N_0(D_G)$. A term in the sum
$$
\sum_{\gamma_1\gamma_2 \gamma_3 = \gamma}
h(\gamma_1)f(\gamma_2)\overline{h(\gamma_3^{-1})}$$
is zero unless $\gamma_2 = s(\gamma_1)$ and $\gamma =
r(\gamma_1)$. Since $U$ is $G$-invariant this shows that $hfh^* \in C_0(U)$.
\end{proof}
\end{cor}

In comparison with the condition
  for simplicity which can be derived from Renault's work, note that
  although the statement does not appear explicitly in \cite{Re1} his
  methods can give the conclusion in Corollary \ref{invopne}, that
  simplicity is equivalent to the absence of any non-trivial open
  $G$-invariant subset in $G^0$, under the
  assumption that points with trivial isotropy is dense in $G^0$. So
  what we do in Corollary \ref{invopne} is to reduce the assumption,
  and in fact to the least possible. Any discrete group whose reduced group $C^*$-algebra is not simple is
an example which shows that in general the existence of at least one unit
with trivial isotropy can not be omitted in Theorem
\ref{simplicity1}.

In a weak moment one might hope that there is a bijection between the ideals of $C^*_r(G)$
and the $G$-invariant ideals of $D_G$ in the setting
of Theorem \ref{simplicity1}, but elementary examples such as the
product of a
discrete group and a locally compact Hausdorff space, shows that
this is certainly not the case. Theorem \ref{simplicity1} is only a
result on the presence or absence of ideals in $C^*_r(G)$.

\subsection{Discrete abelian isotropy and Cartan subalgebras}

In general an abelian $C^*$-subalgebra $D$ of a given $C^*$-algebra
$A$ is \emph{regular} when $A$ is generated as a $C^*$-algebra
by $N(D)$. Following Renault, cf. \cite{Re3}, we say that $D$ is a \emph{Cartan
  subalgebra} in $A$ when
\begin{enumerate}
\item[(i)] $D$ contains an approximate unit in $A$;
\item[(ii)] $D$ is maximal abelian;
\item[(iii)] $D$ is regular, and
\item[(iv)] there exists a faithful conditional expectation $Q : A \to
  D$ of $A$ onto $D$.
\end{enumerate}

Returning to the case where $A = C^*_r(G)$ and 
$D = D_G = C^*_r(G) \cap
  B_0(G^0)$, it follows from Lemma \ref{condexp} that $P_G$ is a faithful conditional
expectation of $C^*_r(G)$ onto $D_G$, and from Lemma \ref{type1} that every $n
\in \alg^* G$ which is a supported in a bisection is a
$D_G$-normalizer. This shows that $D_G$ is regular. It is easy to see
that $D_G$ contains an approximate unit for $C^*_r(G)$, cf. the proof
of Theorem \ref{cartan}, and there is therefore only one thing
missing in Renault's definition of a Cartan subalgebra from \cite{Re3}:
In general $D_G$ is not maximal abelian. In this section we
impose additional conditions on $G$ which hold in many of the applications of
the theory to dynamical systems and which
ensure that $D_G$ is a subalgebra of a larger abelian $C^*$-algebra
which \emph{is} a Cartan subalgebra in the sense of Renault.

Set
$$
\Is G = \left\{\gamma \in G: \ r(\gamma) = s(\gamma) \right\}
$$
which is sometimes called \emph{the isotropy bundle} of $G$. Note that
$\Is G$ is a closed sub-groupoid of $G$.  In the
following we often assume that $\Is(G)$ is abelian, i.e. that
$\gamma_1\gamma_2 = \gamma_2\gamma_1$ for all $G^{(2)} \cap \left(\Is G
  \times \Is G\right)$.

Set
$$
D'_G = \left\{ a \in C^*_r(G): \ \supp j(a) \subseteq \Is(G)
\right\} . 
$$

\begin{lemma}\label{DIsformula} $D'_G$ is a $C^*$-subalgebra of
  $C^*_r(G)$. In fact,
\begin{equation}\label{commutant}
D'_G = \left\{ a \in C^*_r(G): \ ah = ha \ \forall h \in C_c(G^0)
\right\} .
\end{equation}
\begin{proof} It suffices to prove (\ref{commutant}). Let $h \in C_c(G^0)$, $a \in C^*_r(G)$. Then
  $j(ah)(\gamma) = j(a)(\gamma)h\left(s(\gamma)\right)$ and
  $j(ha)(\gamma) = h\left(r(\gamma)\right)j(a)(\gamma)$ for all
  $\gamma \in G$ by Lemma \ref{onj}. Hence $j(ah) = j(ha)$ when $a \in
  D'_G$ and by Corollary \ref{jinj} this implies that $ah = ha$.

Assume next that $a \in C^*_r(G)$ commutes with every element of
  $C_c(G^0)$ and consider an element $\gamma \in G$ with $j(a)(\gamma)
  \neq 0$. If $\gamma \notin \Is G$ we can pick an element $h \in C_c(G^0)$ such that $h(r(\gamma)) = 0$
while $h(s(\gamma)) = 1$. Then $j(ha)(\gamma) =
h(r(\gamma))j(a)(\gamma) = 0$ while $j(ah)(\gamma) =
j(a)(\gamma)h(s(\gamma)) = j(a)(\gamma) \neq 0$, proving that $j(ah -
ha) \neq 0$. By Corollary \ref{jinj} this implies that $ah \neq
ha$, contradicting our assumption on $a$. It follows that
$j(a)(\gamma) = 0$ for all $\gamma \in G \backslash \Is G$, i.e. $a
\in D'_G$.

\end{proof}
\end{lemma}

\begin{defn} We say that $\Is G$ is \emph{discrete} when $\Is G
  \backslash G^0$ is discrete in the topology inherited from $G$.
\end{defn}

\begin{lemma}\label{condexp2} Assume that $\Is G$ is discrete. It follows that there
  is a faithful surjective conditional expectation $Q : C^*_r(G) \to D'_G$.
\begin{proof} Set $\Is_{ess} G = \left\{ \gamma \in \Is G : \
    j(a)(\gamma) \neq 0 \ \text{for some} \ a \in C^*_r(G) \right\}
$. Let $\gamma \in \Is_{ess} G \backslash G^0$. Since $\Is G$ is discrete in the
  topology inherited from $G$ there is a bisection $U$ such that
  $U \cap \Is G = \left\{\gamma\right\}$. Since $\gamma \in \Is_{ess} G$ there
  is also an element $a \in C^*_r(G)$ such that $j(a)(\gamma) = 1$. Let $h
  \in C_c(G)$ be supported in $U$ such that $h(\gamma) =1$. By
  Lemma \ref{product2} there is an element $a_{\gamma} = h \cdot a \in
  C^*_r(G)$ such that
  $j(a_{\gamma}) = 1_{\left\{\gamma\right\}}$. It follows from Corollary \ref{jinj}
  that $a_{\gamma} = 1_{\{\gamma\}}$. I.e. we have shown that
  $1_{\{\gamma\}} \in D'_G$ when $\gamma \in \Is_{ess} G \backslash
  G^0$. 

When $f \in \alg^* G$
the set $\supp f \cap \left(\Is G \backslash G^0\right)$ is finite and
we set
\begin{equation}\label{Qformula}
Q(f) \  = \ f|_{\Is G} \ = \ f|_{G^0} \ \ + \sum_{\gamma \in \Is_{ess} G \backslash G^0} f(\gamma)
1_{\{\gamma\}} .
\end{equation}
Then $Q(f) \in D'_G$. To estimate the norm of $Q(f)$ in $D'_G$,
observe that for every $x \in G^0$ we have a direct sum decomposition
$$
l^2\left(s^{-1}(x)\right) = \oplus_{y \in G^0} \ l^2\left(s^{-1}(x) \cap
  r^{-1}(y)\right)
$$
which is respected by $\pi_x(g)$ when $g \in
B_c\left(\Is G \right)$. It follows that
\begin{equation}\label{oz1}
\left\|\pi_x(g)\right\| = \sup_{y \in G^0} \left\|
  \pi_x(g)|_{l^2\left(s^{-1}(x) \cap r^{-1}(y)\right)}\right\| .
\end{equation}   
Consider a $y \in G^0$ such that $r^{-1}(y) \cap
s^{-1}(x) \neq \emptyset$ and choose $\gamma_0 \in r^{-1}(y) \cap
s^{-1}(x)$. We define a unitary $V : l^2\left(r^{-1}(y) \cap
  s^{-1}(x)\right) \to l^2\left(r^{-1}(y) \cap \Is G\right)$ such that
$V\psi(\eta) = \psi(\eta \gamma_0)$. Then 
$$
V\pi_x(g)V^* = \pi_{y}(g)|_{l^2(r^{-1}(y) \cap \Is G)}
$$
and hence (\ref{oz1}) implies that 
$\left\|\pi_x(g)\right\| \leq \sup_{y \in G^0}  \left\|
    \pi_{y}(g)|_{l^2(r^{-1}(y) \cap \Is G)}\right\|$. 
It follows first that
$\left\|\pi_x(g)\right\| \leq \left\|g\right\|_{B^*_r\left(\Is G\right)}$,
and then that $\left\|g\right\|_{B^*_r(G)} \leq \left\|g
\right\|_{B^*_r(\Is G)}$. Since the reversed inequality is obvious, we
conclude that $\left\|g\right\|_{B^*_r(G)} = \left\|g
\right\|_{B^*_r(\Is G)}$. In particular, $\left\|Q(f)\right\|_{D'_G} = \left\|Q(f)
\right\|_{B^*_r(\Is G)}$. 

Let $y \in G^0$ and note that
$$
\left\|\pi_y\left(Q(f)\right)\right\|_{l^2\left(s^{-1}(y) \cap \Is G\right)} = \left\|
  P_y \pi_y(f)P_y\right\|_{l^2(s^{-1}(y))},
$$
where $P_y : l^2\left(s^{-1}(y) \right) \to l^2\left(s^{-1}(y) \cap
  \Is G\right)$ is the orthogonal projection. It follows that
$$
\left\|\pi_y\left(Q(f)\right)\right\|_{l^2\left(s^{-1}(y) \cap \Is G\right)} \leq
\left\|\pi_y(f)\right\|_{l^2(s^{-1}(y))}.
$$ 
Since $y \in G^0$ was arbitrary we conclude that
$$
\left\|Q(f)\right\|_{D'_G} = \left\|Q(f)
\right\|_{B^*_r(\Is G)} \leq \left\|f\right\|_{C^*_r (G)} .
$$
 Hence $Q$ extends by continuity to a linear map $Q: C^*_r(G) \to
  D'_G$ of norm $1$.

Let $a \in D'_G$. Choose a sequence $\{f_n\} \subseteq \alg^* G$
such that $\lim_{n \to \infty} f_n = a$ in $C^*_r(G)$. Then
$\lim_{n \to \infty} Q(f_n) = Q(a)$. Furthermore,
$$
j(Q(a))(\gamma) = \lim_{n\to \infty} Q(f_n)(\gamma) = \begin{cases} 0,
& \text{when} \ \gamma \notin \Is_{ess} G \\ \lim_{n\to \infty}
f_n(\gamma) = j(a)(\gamma), & \ \text{when} \ \gamma \in \Is_{ess} G. 
\end{cases}
$$
This shows that $j(Q(a)) = j(a)$, and it follows then from Corollary
\ref{jinj} that $Q(a) = a$. Thus $Q :
C^*_r(G) \to D'_G$ is a linear surjective idempotent map of norm
one. It is easy to check that $Q$ is also positive and hence
a conditional expectation. $Q$ is faithful because $P_G \circ Q = P_G$ and
$P_G$ is faithful by Lemma \ref{condexp}.
\end{proof}
\end{lemma}

\begin{cor}\label{acorollary} Assume that $\Is G$ is discrete. Then
$$
D'_G = \overline{C^*_r(G) \cap
    B_c\left(\Is G\right)}.
$$
\begin{proof} The inclusion $C^*_r(G) \cap B_c(\Is G) \subseteq
  D'_G$ is obvious and it follows from Lemma \ref{condexp2} and
  (\ref{Qformula}) that $C^*_r(G) \cap B_c(\Is G)$ is dense in $D'_G$.
\end{proof}
\end{cor}

\begin{thm}\label{cartan}  Assume that $\Is G$ is abelian and
  discrete.  It follows that
  $D'_G$ is a Cartan subalgebra of $C^*_r(G)$.
\begin{proof}  Let $a,b \in D'_G$. Since $\Is G$ is abelian it
  follows from Lemma \ref{onj} that $j(ab) = j(ba)$. By Corollary
  \ref{jinj} this implies that $ab = ba$, proving that $D'_G$ is
  abelian. We check the conditions (i) through (iv) which were listed
  at the beginning of this section: To check condition (i), note that
$C_c\left(G^0\right) \subseteq D'_G$
by Lemma \ref{rdiscrete}.  It is elementary to check that a bounded and
increasing net of non-negative functions from $C_c(G^0)$ which
eventually become constant $1$ on every compact subset of $G^0$ will be an approximate
unit relative to elements from $\alg^* G$ and hence on
all of $C^*_r(G)$.

(ii) follows from (\ref{commutant}). 

To establish (iii) it suffices to show that an element $f \in C_c(G)$
which is supported in a bisection is a $D'_G$-normalizer. Let $a \in D'_G$, $\gamma \in G$. By Lemma \ref{onj}
$$
j(f^*af)(\gamma) = \sum_{\gamma_1 \gamma_2 \gamma_3 =
  \gamma}
\overline{f\left(\gamma_1^{-1}\right)}j(a)\left(\gamma_2\right)f\left(\gamma_3\right)
.
$$
Since $j(a)$ is supported in $\Is G$, $\overline{f\left(\gamma_1^{-1}\right)}j(a)\left(\gamma_2\right)f\left(\gamma_3\right)$ is
zero unless $r\left(\gamma_3\right) = s\left(\gamma_2\right) =
r\left(\gamma_2\right) = s\left(\gamma_1\right)$. Since
$r$ and $s$ are both injective on $\supp f$ there is an (injective)
map $\theta : r(\supp f)
\to s(\supp f)$ such that $f(\mu) = 0$ unless $\theta\left(r(\mu)\right) =
s(\mu)$. So if $\overline{f\left(\gamma_1^{-1}\right)}j(a)\left(\gamma_2\right)f\left(\gamma_3\right)$ is not zero we must also
have that $\theta\left(r(\gamma_3)\right) = s\left(\gamma_3\right)$
and $\theta\left(s\left(\gamma_1\right)\right) =
\theta\left(r\left(\gamma_1^{-1}\right)\right) =
s\left(\gamma_1^{-1}\right) = r\left(\gamma_1\right)$. As observed we
must also have that $s\left(\gamma_1\right) = r\left(\gamma_3\right)$
and it follows that $s(\gamma_3) =r\left(\gamma_1\right)$. Since  $r(\gamma) =
r\left(\gamma_1\right)$  and $s\left(\gamma_3\right) = s(\gamma)$ this
implies that $r(\gamma) =
s(\gamma)$. Thus $j(f^*af)$ is supported in $\Is G$, i.e. $f^*af \in D'_G$.

(iv) follows from Lemma \ref{condexp2}. 
 
\end{proof}
\end{thm}  

A semi \'etale groupoid $G$ is a \emph{semi \'etale equivalence
  relation} when $\Is G = G^0$. To distinguish these groupoids from
the more general ones we shall denote a semi \'etale equivalence
relation by $R$.

\begin{lemma}\label{relation} Let $R$ be a semi \'etale equivalence
  relation and $P_R : C^*_r(R) \to D_R$ the corresponding conditional expectation. Let $\epsilon >
  0$ and $a \in C^*_r(R)$ be given. It follows that there are positive
  elements $d_i, i = 1,2, \dots, N$, in $C_c\left(R^0\right)
  \subseteq D_R$ such that
$$
\left\| P_R(a) - \sum_{i=1}^N d_i a d_i\right\| \leq \epsilon .
$$
\begin{proof} Choose $ f \in \alg^* R$ such that $\left\|a - f\right\|
  \leq {\epsilon}$. Set $E = R \backslash R^0$. Since there is
  only
  trivial isotropy in $R$ we can cover $E$ by open bisections $U$ such
  that $r(U) \cap s(U) = \emptyset$. We can therefore apply Lemma
  \ref{techlemma} to obtain a decomposition $f = f|_{R^0} \plus
  \sum_{j=1}^M h_j$ where $h_j \in \alg^* R$ is supported in a
  bisection $U_j$ with $r\left(U_j\right) \cap s\left(U_j\right) =
  \emptyset$. Set $K = \left(\supp f|_{R^0}\right) \cup
  \bigcup_{j=1}^M  \left(r\left(\supp h_j \right) \cup
    s\left(\supp h_j\right)\right)$ and cover $K$ by a finite open cover
  $V_i, i = 1,2, \dots, N$, such that $V_i \cap r\left(\supp h_j\right) \neq
  \emptyset \ \Rightarrow \ V_i \cap s\left(\supp h_j\right) = \emptyset$
  for all $i,j$, and let $k_i \in C_c\left(R^0\right), i = 1,2, \dots, N$, be a partition of
  unity on $K$ subordinate to $V_i, i=1,2, \dots, N$. Set $d_i =
  \sqrt{k_i}$. Then $P_R(f) = f|_{R^0} = \sum_{i=1}^N d_ifd_i$ and
  $\left\| P_R(a) - \sum_{i=1}^N d_iad_i\right\| \leq 2 \epsilon$. 
 \end{proof}
\end{lemma}


\section{The $C^*$-algebra of a locally injective map}

Let $X$ and $Y$ be compact metrizable Hausdorff spaces and $\varphi :
X \to Y$ a continuous and
locally injective map. Set
$$
R(\varphi) = \left\{ (x,y) \in X \times X : \ \varphi(x) = \varphi(y)
\right\} .
$$
This is a semi \'etale relation in the topology inherited from $X
\times X$ and we present the information on the structure of
$C^*_r\left(R(\varphi)\right)$ which will be needed in the subsequent sections.

\begin{lemma}\label{normest} Set $K = \max_x
\#  \varphi^{-1}\left(\varphi(x)\right)$. Then
$$
\left\|d\right\| \leq K \sup_{(x,y) \in R(\varphi)}
\left|j(d)(x,y)\right| 
$$
for all $d \in B^*_r\left(R(\varphi)\right)$.
\begin{proof} Let $\{f_n\} \subseteq B_c\left(R(\varphi)\right)$ be a
  sequence converging to $d$ in $B^*_r\left(R(\varphi)\right)$. (The
  subscript $c$ is redundant in this case because $R(\varphi)$ is
  compact, but we keep it on for consistency with the notation of
  Section \ref{Sec1}.) Then
  $\{f_n\}$ also converges to $j(d)$, uniformly on $R(\varphi)$, by
  (\ref{on1}). It
  suffices therefore to prove the desired inequality when $d \in
  B_c\left(R(\varphi)\right)$: When $x \in X$ and
  $\psi \in l^2\left(s^{-1}(x)\right)$ we find that
\begin{equation*}
\begin{split}
&\left\|\pi_x(d)\psi\right\|^2 =  \sum_{y \in s^{-1}(x)} \left|\sum_z
  d(y,z)\psi(z)\right|^2 \leq \sum_{y \in s^{-1}(x)} \|\psi\|^2_{l^2(s^{-1}(x))} \sum_z 
  \left|d(y,z)\right|^2 \\
&\leq K^2 \|\psi\|^2_{l^2(s^{-1}(x))}
\left(\sup_{(x,y) \in R(\varphi)} \left|d(x,y)\right|\right)^2 .
\end{split}
\end{equation*}
\end{proof}
\end{lemma}

We will consider the elements of $B_r^*(R(\varphi))$ as functions on
$R(\varphi)$, as we can by Corollary \ref{jinj}. It follows then from Lemma
\ref{normest} that the $C^*$-norm of $B^*_r(R(\varphi))$ is equivalent
to the supremum norm.

Set
$$
Y \times_{\varphi} X = \left\{ (a,x) \in Y \times X : \ a = \varphi(x) \right\} 
$$
which is a closed subset of $\varphi(X) \times X$. Let
$B(Y\times_{\varphi} X)$ denote the vector space of bounded complex
functions on $Y \times_{\varphi} X$. We intend to construct an
imprimitivity bimodule, in the sense of Rieffel, out of $B\left(Y
  \times_{\varphi} X\right)$, and we refer to \cite{RW} for a nice
exposition of the theory we rely on.

When $h,k \in B\left(Y \times_{\varphi} X\right)$ we
define 
$\left< h, k \right> : \ R(\varphi) \to \mathbb C$
such that 
$$
\left< h,k\right> (x,y) =
\overline{h\left(\varphi(x),x\right)}k\left(\varphi(y),y\right) ,
$$
and
$(h,k) : \varphi(X) \to \mathbb C$
such that 
$$
(h,k)(a) = \sum_{z \in \varphi^{-1}(a)} h(a,z)\overline{k(a,z)} .
$$
Note that $\left<h,h\right>$ is positive in
$B^*_r\left(R(\varphi)\right)$ since
$\pi_z\left(\left<h,h\right>\right)$ is positive as an operator on
$l^2\left(s^{-1}(z)\right)$ for every $z \in X$.

Let $B(\varphi(X))$ be the $C^*$-algebra of bounded complex functions
on $\varphi(X)$. When $g \in B(\varphi(X))$ and $h \in  B\left(Y \times_{\varphi}
  X\right)$ we define $gh\in B\left(Y \times_{\varphi}
  X\right)$ such that
$$
gh(a,x) = g(a)h(a,x),
$$ 
and when $h \in  B\left(Y \times_{\varphi}
  X\right)$, $f \in B_c\left(R(\varphi)\right)$ we define $hf \in  B\left(Y \times_{\varphi}
  X\right)$ such that
$$
hf(a,x) =\sum_{z \in \varphi^{-1}(a)}h(a,z)f(z,x) .
$$
It is straightforward to check that with these definitions $B\left(Y
  \times_{\varphi} X\right)$ has all
the properties of a
  $B(\varphi(X))$-$B^*_r\left(R(\varphi)\right)$-imprimitivity bimodule,
  except that the fullness of the Hilbert modules may fail. See
  Definition 3.1 on page 42 of \cite{RW}.

Let
$
E_{\varphi}$
be the Hilbert $C^*_r\left(R(\varphi)\right)$-module   
$$
E_{\varphi} = \overline{\Span \left\{ fg : \ f \in C\left(Y
    \times_{\varphi} X\right), \ g \in C^*_r\left(R(\varphi)\right)
\right\}} .
$$
In particular, $\left<h,k\right> \in C^*_r\left(R(\varphi)\right)$ when $h,k \in
E_{\varphi}$. Set
$$
A_{\varphi} = \overline{\Span \left\{ (h,k) : \ h,k \in E_{\varphi}
  \right\}} ,
$$
which is unital $C^*$-subalgebra of $B(\varphi(X))$. Note that
$E_{\varphi}$ is then a
full left Hilbert $A_{\varphi}$-module.

\begin{thm}\label{morthm2} $E_{\varphi}$ is a full Hilbert
  $C^*_r\left(R(\varphi)\right)$-module and hence an
  $A_{\varphi}$-$C^*_r\left(R(\varphi)\right)$-imprimitivity
  bimodule.
\begin{proof} We must show that the closed twosided ideal of
  $C^*_r\left(R(\varphi)\right)$ generated by
\begin{equation}\label{AV17}
\left\{ \left< h,k\right> : \ h,k \in C\left(Y \times_{\varphi}
    X\right) \right\}
\end{equation}
is all of $C^*_r\left(R(\varphi)\right)$. By definition
$C^*_r\left(R(\varphi)\right)$ is generated as a $C^*$-algebra by
$C\left(R(\varphi)\right)$, and it suffices therefore to show that
$C(R(\varphi))$ is contained in the closed linear span of the elements
from (\ref{AV17}). Let $f \in C\left(R(\varphi)\right)$. By
Tietze's extension theorem there is a $g \in C(X \times X)$ such that
$g|_{R(\varphi)} = f$ and we can therefore approximate $f$ in the
supremum norm, and hence also in the $C^*$-norm of
$C^*_r\left(R(\varphi)\right)$ by a linear combination of functions
of the form $\mu \otimes \kappa$, where $\mu,\kappa \in C(X)$ and
$\mu \otimes \kappa (x,y) = {\mu(x)}\kappa(y)$.
Define $h,k : Y\times_{\varphi} X \to \mathbb C$ such that $h(a,x) =
\overline{\mu(x)}$ and $k(a,x) = \kappa(x)$. Since $\left<h,k\right> =
\mu\otimes \kappa$, this completes the proof.
  \end{proof}
\end{thm}

\begin{cor}\label{Hcor}
Let $Z_{\varphi}$ denote the Gelfand spectrum of $A_{\varphi}$. It
follows that there is an $n \in \mathbb N$ and a projection $p \in
C\left(Z_{\varphi},M_n(\mathbb C)\right)$ such that
$$
C^*_r\left(R(\varphi)\right) \simeq pC\left(Z_{\varphi},M_n(\mathbb C)\right)p .
$$
\begin{proof} It follows from Theorem \ref{morthm2} that
  $C^*_r\left(R(\varphi)\right)$ is Morita equivalent to
  $C\left(Z_{\varphi}\right)$ and then from \cite{BGR} and \cite{Br} that
  $C^*_r\left(R(\varphi)\right)$ is a corner in $C(Z_{\varphi})
  \otimes \mathbb K$, where $\mathbb K$ denotes the $C^*$-algebra of
  compact operators on an infinite dimensional separable Hilbert space. Such a corner has the form
  $pC\left(Z_{\varphi},M_n(\mathbb C)\right)p$ for some $n$ and some
  $p \in C\left(Z_{\varphi},M_n(\mathbb C)\right)$.
\end{proof}
\end{cor}

In the terminology from the classification program of $C^*$-algebras,
cf. e.g. \cite{EGL},
what Corollary \ref{Hcor} says is that $C^*_r\left(R(\varphi)\right)$ is a
direct sum of homogeneous $C^*$-algebras. In the context of type I
$C^*$-algebras, cf. e.g. \cite{RW}, it says that
$C^*_r\left(R(\varphi)\right)$ is a direct sum of $n$-homogeneous
$C^*$-algebras with trivial Dixmier-Douady invariant. In
particular the primitive ideal space of $C^*_r\left(R(\varphi)\right)$
is Hausdorff. It is possible, but lengthy to give a complete description of the
primitive ideal space. Only the following partial results in that direction will
be needed.

Let $a \in \varphi(X)$. We can then define a $*$-homomorphism
$\tilde{\psi}_a : B_c\left(R(\varphi)\right) \to M_{\varphi^{-1}(a)}(\mathbb C)$
such that $\tilde{\psi}_a(h) = h|_{\varphi^{-1}(a) \times
  \varphi^{-1}(a)}$. Since $\left\|h|_{\varphi^{-1}(a) \times
  \varphi^{-1}(a)}\right\|_{ M_{\varphi^{-1}(a)}(\mathbb C)} =
\left\|\pi_x(h)\right\|$ for any $x \in \varphi^{-1}(a)$ we find that 
$\left\|h|_{\varphi^{-1}(a) \times
  \varphi^{-1}(a)}\right\|_{M_{\varphi^{-1}(a)}(\mathbb C)} \leq
\left\|h\right\|$
for all $h \in B_c\left(R(\varphi)\right)$. Hence $\tilde{\psi}_a$ extends to a
$*$-homomorphism $\tilde{\psi}_a : B^*_r\left(R(\varphi)\right) \to
M_{\varphi^{-1}(a)}(\mathbb C)$ which is clearly surjective. Set
$\psi_a = \tilde{\psi}_a|_{C^*_r\left(R(\varphi)\right)}$ which is
also surjective since $C\left(R(\varphi)\right) \subseteq
C^*_r\left(R(\varphi)\right)$.

\begin{lemma}\label{iotadense} $\left\{ \ker \psi_a : \ a \in
    \varphi(X) \right\}$ is dense in the primitive ideal space $\Prim
   C^*_r\left(R(\varphi)\right)$ of $C^*_r\left(R(\varphi)\right)$.
\begin{proof} Let $W$ be a non-empty open subset of $\Prim
  C^*_r\left(R(\varphi)\right)$. Since $\Prim
  C^*_r\left(R(\varphi)\right)$ is Hausdorff there is an element $d \in
  C^*_r\left(R(\varphi)\right)$ such that $d \neq 0$ and
$$
\left\{ \ker \pi \in \Prim
  C^*_r\left(R(\varphi)\right) : \ \pi(d) \neq 0 \right\} \subseteq W,
$$
cf. e.g. \cite{RW}. Since $d \neq 0$ there is an $a \in \varphi(X)$ such that
$d|_{\varphi^{-1}(a) \times \varphi^{-1}(a)} \neq 0$. It follows that $\ker
\psi_a \in W$.
\end{proof}
\end{lemma}

\begin{lemma}\label{centralproj} For each $j = 1,2,3, \dots $ there is
  a (possibly zero) projection $p_j$ in the center of
  $C^*_r\left(R(\varphi)\right)$ such that
$$
p_j(x,y) = \begin{cases} 1, & \ \text{when $x=y$ and $\#
    \varphi^{-1}(\varphi(x)) = j$}, \\ 0, &
    \text{otherwise}. \end{cases}
$$
\begin{proof} Let $p$ be the function on $R(\varphi)$ which is constant
  $1$. Then $p \in C\left(R(\varphi)\right) \subseteq
  C^*_r\left(R(\varphi)\right)$ and hence $p \star p^*|_X \in
  C^*_r\left(R(\varphi)\right)$. Since 
$$
p\star p^*(x,x) = \sum_{j=1}^{\infty} jp_j(x,x),
$$
it follows from spectral theory that $p_j \in
C^*_r\left(R(\varphi)\right)$ for all $j$. It is straigthforward to
check that $p_j$ is central.
\end{proof}
\end{lemma}

Note that there are only finitely many $j \in \mathbb N$ for which
$p_j \neq 0$ and that $\sum_j p_j =1$.

\begin{lemma}\label{cruxlemma} Let $z \in \varphi(X)$ and set $j = \#
  \varphi^{-1}(z)$. There is then an
  open neighborhood $U$ of $z$ and open sets $V_i, i=1,2, \dots,j$, in
  $X$, such
  that
\begin{enumerate}
\item[1)] $\varphi^{-1}\left(\overline{U}\right) \subseteq V_1 \cup V_2
  \cup \dots \cup V_j$,
\item[2)] $\overline{V_i} \cap \overline{V_{i'}} = \emptyset, \ i \neq
  i'$, and
\item[3)] $\varphi$ is injective on $\overline{V_i}$ for
  each $i$.
\end{enumerate} 
\begin{proof} Since $\varphi$ is locally injective there are open sets
  $V_i, i = 1,2, \dots, j$, such that 2) and 3) hold and
\begin{equation}\label{contraeq}  
\varphi^{-1}(z) \subseteq \bigcup_{i=1}^j V_i. 
\end{equation}
If there is no open
  neighborhood $U$ of $z$ for which 1) holds there is a sequence
  $\left\{x_{n}\right\} \subseteq X \backslash \bigcup_{i=1}^j
  V_i$ such that $\lim_{n} \varphi\left(x_{n}\right) = z$. A
  condensation point $x$ of this sequence gives us an element $x \in X \backslash \bigcup_{i=1}^j
  V_i$ such that $\varphi(x) = z$, contradicting (\ref{contraeq}).  
\end{proof}
\end{lemma}

For each $j \in \mathbb N$, set 
$$
L_j = \left\{(x,y) \in R(\varphi) : \ \# \varphi^{-1}(\varphi(x)) = j
\right\} .
$$

\begin{lemma}\label{contonlj} Let $a \in
  C^*_r\left(R(\varphi)\right)$. Then $a|_{L_j}$ is continuous on
  $L_j$ for every $j$.
\begin{proof} Since continuity is preserved under uniform convergence
  it suffices to prove this when $a \in \alg^* R(\varphi)$, and
  hence in fact when $a = f_1 \star f_2 \star \dots \star f_N$ for
  some $f_1,f_2, \dots, f_N \in C\left(R(\varphi)\right)$. Let $(x,y)
  \in L_j$ and set $z = \varphi(x) = \varphi(y)$. Let $U$ and $V_i, i
  = 1,2,\dots, j$, be as in Lemma \ref{cruxlemma}.
For every $z' \in U$ with $\# \varphi^{-1}\left(z'\right) =
j$ there are unique elements $\lambda_k(z') \in V_k$ such that
$\varphi^{-1}(z') = \left\{ \lambda_1(z'),\lambda_2(z'), \dots,
  \lambda_j(z')\right\}$. Then
\begin{equation*}
\begin{split}
&f_1 \star f_2 \star \dots \star f_N(x',y') = \\
&\sum_{k_1,k_2, \dots,
  k_{N-1} }
f_1\left(x',\lambda_{k_1}(\varphi(x'))\right)f_2\left(\lambda_{k_1}\left(\varphi(x')\right),\lambda_{k_2}\left(\varphi(x')\right)\right)
\cdots\\
& \ \ \ \ \ \ \ \ \ \ \ \ \ \ \ \ \ \ \ \ \ \ \ \ \cdots
f_{N-1}\left(\lambda_{k_{N-2}}\left(\varphi(x')\right),\lambda_{k_{N-1}}\left(\varphi(x')\right)\right)
f_{N}\left(\lambda_{k_{N-1}}\left(\varphi(x')\right),y'\right)
\end{split}
\end{equation*}
when $(x',y') \in L_j$ and $x' \in \varphi^{-1}(U)$. It suffices therefore
to prove that each $\lambda_k$ is continuous on $U \cap \left\{ z \in
  Y: \ \# \varphi^{-1}(z) = j\right\}$. Let $\{a_n\}$ be a
sequence in $U \cap \left\{ z \in
  Y: \ \# \varphi^{-1}(z) = j\right\}$ converging to $a \in U \cap \left\{ z \in
  Y: \ \# \varphi^{-1}(z) = j\right\}$. If
$\left\{\lambda_k(a_n)\right\}$ does not converge to $\lambda_k(a)$ for
some $k$, the sequence $\left\{\lambda_k(a_n)\right\}$ will have a condensation point $x \in
\overline{V_k} \backslash \left\{\lambda_k(a)\right\}$. Since $\varphi(x) =
\varphi\left(\lambda_k(a)\right) = a$, this contradicts condition 3)
of Lemma \ref{cruxlemma}.
\end{proof}
\end{lemma}

\section{Dynamical systems}\label{sec2}

Let $X$ be a compact metrizable Hausdorff space and $\varphi : X \to X$
a continuous map. We assume that $\varphi$ is locally injective. Set
$$
\Gamma_{\varphi} = \left\{ (x,k,y) \in X \times \mathbb Z \times X : \ \exists
  a,b \in \mathbb N, \ k = a-b, \ \varphi^a(x) = \varphi^b(y) \right\} .
$$
This is a groupoid with the set of composable pairs being
$$
\Gamma_{\varphi}^{(2)} \ =  \ \left\{\left((x,k,y), (x',k',y')\right) \in \Gamma_{\varphi} \times
  \Gamma_{\varphi} : \ y = x'\right\}.
$$
The multiplication and inversion are given by 
$$
(x,k,y)(y,k',y') = (x,k+k',y') \ \text{and}  \ (x,k,y)^{-1} = (y,-k,x)
.
$$
To turn $\Gamma_{\varphi}$ into a locally compact topological groupoid, fix $k \in \mathbb Z$. For each $n \in \mathbb N$ such that
$n+k \geq 0$, set
$$
{\Gamma_{\varphi}}(k,n) = \left\{ \left(x,l, y\right) \in X \times \mathbb
  Z \times X: \ l =k, \ \varphi^{k+i}(x) = \varphi^i(y), \ i \geq n \right\} .
$$
This is a closed subset of the topological product $X \times \mathbb Z
\times X$ and hence a locally compact Hausdorff space in the relative
topology.
Since $\varphi$ is locally injective $\Gamma_{\varphi}(k,n)$ is an open subset of
$\Gamma_{\varphi}(k,n+1)$ and hence the union
$$
{\Gamma_{\varphi}}(k) = \bigcup_{n \geq -k} {\Gamma_{\varphi}}(k,n) 
$$
is a locally compact Hausdorff space in the inductive limit topology. The disjoint union
$$
\Gamma_{\varphi} = \bigcup_{k \in \mathbb Z} {\Gamma_{\varphi}}(k)
$$
is then a locally compact Hausdorff space in the topology where each
${\Gamma_{\varphi}}(k)$ is an open and closed set. In fact, as is easily verified, $\Gamma_{\varphi}$ is a locally
compact groupoid in the sense of \cite{Re1}. In the following we shall often identify the unit space $\Gamma_{\varphi}^0$ of
$\Gamma_{\varphi}$ with $X$ via the map $x \to (x,0,x)$ which is a
homeomorphism. The local
injectivity of $\varphi$ ensures that the range map $r(x,k,y) = x$ is
locally injective, i.e. $\Gamma_{\varphi}$
is semi \'etale. Note that every isotropy group of $\Gamma_{\varphi}$ is a
subgroup of $\mathbb Z$. In particular, $\Is \Gamma_{\varphi}$ is abelian.

When $\varphi$ besides being locally injective is also open, and hence
a local homeomorphism, $\Gamma_{\varphi}$ is an \'etale groupoid, which was
introduced in increasing generality in \cite{Re1}, \cite{De},
\cite{A} and \cite{Re2}. However, when $\varphi$ is not open $\Gamma_{\varphi}$ is no longer
\'etale, merely
semi \'etale.

\begin{lemma}\label{discreteisotropy} Assume that $\left\{x \in X : \
  \varphi^k(x) = x\right\}$ is discrete in the topology inherited from
$X$ for all $k \in \mathbb
N$. It follows that $\Is \Gamma_{\varphi}$ is discrete.
\begin{proof} Let $\gamma = (x_0,k,x_0) \in \Is \Gamma_{\varphi} \backslash
  \Gamma_{\varphi}^0$. Then $k \neq 0$ and $x_0 \in \Gamma_{\varphi}(k,n)$ for some $n \geq
  1$. Note that $\varphi^n(x_0)$ is $|k|$-periodic. By assumption there
  is an open neighborhood $U$ of $x_0$ such that $x_0$ is the only
  element $x$
  of $U$ for which $\varphi^{n}(x)$ is $|k|$-periodic. Then
$$
W = \left\{ (x,k,y) \in \Gamma_{\varphi}(k,n) : \ x,y \in U \right\}
$$
is an open subset of $\Gamma_{\varphi}$ such that $W \cap \Is \Gamma_{\varphi} = \{\gamma\}$.
\end{proof}
\end{lemma}

\begin{thm}\label{cartan13} Assume that $\left\{x \in X : \
  \varphi^k(x) = x\right\}$ is discrete in the topology inherited from
$X$ for all $k \in \mathbb
N$. It follows that $D'_{\Gamma_{\varphi}}$ is a Cartan subalgebra of
$C^*_r(\Gamma_{\varphi})$.
\begin{proof} This is now a consequence of Theorem \ref{cartan}.
\end{proof}
\end{thm}

\begin{lemma}\label{gaugeaction} Let $\mathbb T$ be the unit circle in
  $\mathbb C$. There is a continuous action $\mathbb
  T \ni z \ \mapsto \ \beta_z \in \Aut C^*_r(\Gamma_{\varphi})$ such that
\begin{equation}\label{actionformula}
\beta_z(f)(x,k,y) = z^kf(x,k,y) 
\end{equation}
when $f \in C_c(\Gamma_{\varphi})$ and $(x,k,y) \in \Gamma_{\varphi}$.  
\begin{proof} It is straightforward to check that the formula
  (\ref{actionformula}) defines an automorphism $\beta_z$ of $B_c(\Gamma_{\varphi})$
  such that $\beta_z\left(\alg^* \Gamma_{\varphi}\right) = \alg^*
  \Gamma_{\varphi}$. To see that $\beta_z$ extends by continuity to
  $C^*_r(\Gamma_{\varphi})$, let $x \in X$ and define a unitary $U_z$ on
  $l^2\left(s^{-1}(x)\right)$ such that
$$
U_z\psi(y,k,x) = z^k \psi(y,k,x) . 
$$
Then $\pi_x\left(\beta_z(a)\right) = U_z\pi_x(a)U_z^*$ and
hence $\left\|\pi_x(\beta_z(a))\right\|_{l^2\left(s^{-1}(x)\right)} =
\left\|\pi_x(a)\right\|_{l^2\left(s^{-1}(x)\right)}$. It follows
that $\beta_z$ extends to an automorphism of $C^*_r(\Gamma)$ for each
$z \in \mathbb T$. To check the continuity of $z \mapsto \beta_z(a)$
for each $a \in C^*_r(\Gamma_{\varphi})$ it suffices to check when $a = f \in
C_c(\Gamma_{\varphi})$ is supported in a bisection inside $\Gamma_{\varphi}(k)$ for some $k \in \mathbb
Z$. In this case we have the estimate
\begin{equation*}\label{etim}
\left\|\beta_z(f) - \beta_{z'}(f)\right\| \leq \left|z^k - z'^k\right|
\sup_{\gamma \in G} \left|f(\gamma)\right|
\end{equation*}
by Lemma \ref{renaultineq}. This proves the continuity of $z \mapsto
\beta_z$.

\end{proof}
\end{lemma}

We will refer to the action $\beta$ from Lemma \ref{gaugeaction} as
the \emph{gauge action}. The fixed point algebra of the gauge action
will be denoted by $C^*_r\left(\Gamma_{\varphi}\right)^{\mathbb T}$.

\subsection{ A crossed product description of $C^*_r\left(\Gamma_{\varphi}\right)$}

Note that
$$
\Gamma_{\varphi}(0) = \left\{ (x,k,y) \in \Gamma_{\varphi} : \ k = 0
\right\}
$$
is an open subgroupoid of $\Gamma_{\varphi}$ and hence a semi \'etale
groupoid in itself. We identify $\Gamma_{\varphi}(0)$ with
$$
\left\{ (x,y) \in X \times X : \ \varphi^i(x) = \varphi^i(y) \
  \text{for some} \ i \in \mathbb N\right\}
$$
under the map $(x,y) \mapsto (x,0,y)$. Note that $\Gamma_{\varphi}(0)$
is a semi \'etale equivalence
relation which we 
denote by $R_{\varphi}$ in the following. Similarly,
$\Gamma_{\varphi}(0,n)$ is an open subgroupoid of
$R_{\varphi} \subseteq \Gamma_{\varphi}$ and a semi \'etale equivalence relation
in itself. In fact, $\Gamma_{\varphi}(0,n)$ is isomorphic, as a semi
\'etale equivalence relation, to the semi \'etale equivalence relation
$R\left(\varphi^n\right)$ corresponding to the locally injective map
$\varphi^n$. The isomorphism is given the map
$R\left(\varphi^n\right) \ni (x,y) \mapsto (x,0,y)$, 
and it induces an isomorphism $C^*_r\left(R\left(\varphi^n\right)\right) \simeq
C^*_r\left(\Gamma_{\varphi}(0,n)\right)$. In the following we suppress
these isomorphisms in the notation and identify
$R\left(\varphi^n\right)$ with $\Gamma_{\varphi}(0,n)$ and
$C^*_r\left(R\left(\varphi^n\right)\right)$ with $
C^*_r\left(\Gamma_{\varphi}(0,n)\right)$. Then
$$
R_{\varphi} = \bigcup_{n \in \mathbb N} R\left(\varphi^n\right) .
$$
It follows from Lemma \ref{openincl} that there are isometric
embeddings
$C^*_r\left( R\left(\varphi^n\right)\right) \subseteq C^*_r\left(
  R\left(\varphi^{n+1}\right)\right) \subseteq C^*_r\left( R_{\varphi}\right)$
for all $n$. Since $C_c\left(R_{\varphi}\right) = \bigcup_n
C_c\left(R\left(\varphi^n\right)\right)$ it follows that
\begin{equation}\label{union}
 C^*_r\left(R_{\varphi}\right) = \overline{\bigcup_n
  C^*_r\left(R\left(\varphi^n\right)\right)} .
\end{equation}

Combined with Corollary \ref{Hcor} this shows that
$C^*_r\left(R_{\varphi}\right)$ is an AH-algebra in the terminology
from the classification program for $C^*$-algebras, cf. e.g
\cite{EGL}.

We assume that $\varphi$ is surjective. The aim is to show that there
is then an
endomorphism of $C^*_r\left(R_{\varphi}\right)$ such that
$C^*_r\left(\Gamma_{\varphi}\right)$ is the crossed product of
$C^*_r\left(R_{\varphi}\right)$ by this endomorphism. In the \'etale
case, where $\varphi$ is open, this crossed product decomposition was
established in \cite{A}. Define $m : X \to \mathbb R$ such that
$$
m(x) = \# \left\{ y \in X : \ \varphi(y) = \varphi(x) \right\} .
$$
It follows from Lemma \ref{centralproj} that $m$ is an element of
$D_{R(\varphi)}$ which is central in
$C^*_r\left(R(\varphi)\right)$. Note that $m$ is positive and invertible.

\begin{lemma}\label{nyfulllemma} For each $k \geq 1$ there is a
  $*$-homomorphism $h_k : C^*_r\left(R\left(\varphi^{k}\right)\right) \to
  C^*_r\left(R\left(\varphi^{k+1}\right)\right)$ such that
\begin{equation}\label{hk}
h_k(f)(x,y) = m(x)^{-\frac{1}{2}} m(y)^{-\frac{1}{2}} f\left(\varphi
  (x), \varphi (y)\right).
\end{equation}
\begin{proof} The formula (\ref{hk}) makes sense for all $f \in
  B_c\left(R\left(\varphi^{k}\right)\right)$ and defines a linear map
  $B_c\left(R\left(\varphi^{k}\right)\right) \to
  B_c\left(R\left(\varphi^{k+1}\right)\right)$ which is continuous for the
  supremum norms. When $f \in C\left(R\left(\varphi^{k}\right)\right)$,
 $$
h_k(f) = m^{-1/2} \star \left[f \circ (\varphi \times \varphi)\right] \star
m^{-1/2} 
$$
which is in $C^*_r\left(R\left(\varphi^{k+1}\right)\right)$ since $m$ is by
Lemma \ref{centralproj}. It suffices
therefore to check that $h_k(f^*) = h_k(f)^*$ and $h_k(f \star g) = h_k(f) \star h_k(g)$
  when $f,g \in B_c\left(R\left(\varphi^{k}\right)\right)$. The first
  property is obvious. We check the second:
\begin{equation*}
\begin{split}
& h_k(f)\star h_k(g)(x,y) \\
&\\
&= \sum_{ \left\{ z \in X :  \ \varphi^{k+1}(z) = \varphi^{k+1}(x)\right\}}
m(x )^{-\frac{1}{2}}m(y)^{-\frac{1}{2}} m(z)^{-1}
f\left(\varphi(x),\varphi (z)\right)g\left(\varphi(z),\varphi(y)\right)\\ 
&\\
& = \sum_{\left\{z \in X: \ \varphi(z) = w \right\}} \sum_{\left\{w
    \in X : \ \varphi^k(w) = \varphi^{k+1}(x)\right\}}
 m(x)^{-\frac{1}{2}}m(y)^{-\frac{1}{2}} m(z)^{-1}
f\left(\varphi (x),w\right)g\left(w,\varphi(y)\right) \\
&\\
& = \sum_{\left\{w \in X : \ \varphi^k(w) = \varphi^{k+1}(x)\right\}}
 m(x)^{-\frac{1}{2}}m(y)^{-\frac{1}{2}} 
f\left(\varphi(x),w\right)g\left(w,\varphi(y)\right) \\
& = h_k(f \star g)(x,y) ,
\end{split}
\end{equation*}
where the surjectivity of $\varphi$ was used for the second equality.

\end{proof}
\end{lemma}

\begin{cor}\label{cor1} The function $X \ni x \mapsto m(\varphi^k(x))$
  is in $C^*_r\left(R\left(\varphi^{k+1}\right)\right)$ for $k = 0,1,2,3, \dots$.
\begin{proof} When $h_k$ is the $*$-homomorphism from Lemma
  \ref{nyfulllemma} we have that $m \circ \varphi^k =
  P_{R\left(\varphi^{k+1}\right)}\left(m^{\frac{1}{2}}\star h_k(m \circ \varphi^{k-1}) \star
  m^{\frac{1}{2}}\right)$. In this way the assertion follows from
Lemma \ref{nyfulllemma} by induction.
\end{proof}
\end{cor}

Note that the diagram
\begin{equation}
\begin{xymatrix}{
C^*_r\left(R\left(\varphi^k\right)\right) \ar@{^{(}->}[r]
\ar[d]_-{h_k}  &
C^*_r\left(R\left(\varphi^{k+1}\right)\right) \ar[d]_-{h_{k+1}} \\
C^*_r\left(R\left(\varphi^{k+1}\right)\right) \ar@{^{(}->}[r]
  &
C^*_r\left(R\left(\varphi^{k+2}\right)\right)
}\end{xymatrix}
\end{equation}
commutes for each $k$ so that we obtain a $*$-endomorphism 
$\widehat{\varphi} : C^*_r\left(R_{\varphi}\right) \to
C^*_r\left(R_{\varphi}\right)$
defined such that $\widehat{\varphi}|_{C^*_r\left(R\left(\varphi^k\right)\right)} = h_k$.

We define a function $v : \Gamma_{\varphi} \to \mathbb C$ such that
$$
v(x,k,y) = \begin{cases} m(x)^{-\frac{1}{2}} & \ \text{when} \ k = 1 \
  \text{and} \ y = \varphi(x) \\ 0 & \ \text{otherwise.} \end{cases}
$$
Then $v$ is the product $v = m^{-\frac{1}{2}} 1_{\Gamma_{\varphi}(1,0)}$ 
in $C^*_r\left(\Gamma_{\varphi}\right)$. In particular, $v \in
C^*_r\left(\Gamma_{\varphi}\right)$. By checking on
$C^*_r\left(R\left(\varphi^k\right)\right)$ one finds that
\begin{equation}\label{adv}
vav^* = \widehat{\varphi}(a)
\end{equation}
for all $a \in C^*_r\left(R_{\varphi}\right)$. Similarly a direct
computation shows that $v$ is an isometry, i.e. $v^*v =1$. Unlike the
\'etale case considered in \cite{A}, the inclusion
$v^*C_r^*\left(R_{\varphi}\right)v \subseteq
C^*_r\left(R_{\varphi}\right)$ can fail in the semi \'etale case.

\begin{thm}\label{generation} Assume that $\varphi$ is surjective. It
  follows that $C^*_r\left(\Gamma_{\varphi}\right)$
  is generated, as a $C^*$-algebra, by the isometry $v$ and
  $C^*_r\left(R_{\varphi}\right)$. In fact,
  $C^*_r\left(\Gamma_{\varphi}\right)$ is the crossed product
$$
C^*_r\left(R_{\varphi}\right) \times_{\widehat{\varphi}} \mathbb N
$$
in the sense of Stacey, \cite{St}, and Boyd,
Keswani and Raeburn, \cite{BKR}.
\begin{proof} By definition $C^*_r\left(\Gamma_{\varphi}\right)$ is
  generated by
$$
\bigcup_{n,l \in \mathbb N} C_c\left(\Gamma_{\varphi}(l,n)\right) \cup C_c\left(\Gamma_{\varphi}(-l,n)\right) 
$$
so to prove the first assertion it suffices to show that $C_c\left(\Gamma_{\varphi}(l,n)\right)$ and
$C_c\left(\Gamma_{\varphi}(-l,n)\right)$ are both subsets of the $C^*$-algebra generated by $v$ and
  $C^*_r\left(R_{\varphi}\right)$ for every $l,n$. Assume that $f \in
  C_c\left(\Gamma_{\varphi}(l,n)\right)$. Define the function $g :
  R\left(\varphi^n\right) \to \mathbb C$ such that
$$
g(x,y) =
f\left(x,l,\varphi^l(y)\right)m\left(\varphi^{l-1}(y)\right)^{-\frac{1}{2}}m\left(\varphi^{l-2}(y)\right)^{-\frac{1}{2}}
\cdots m(y)^{-\frac{1}{2}} .
$$
It follows from Corollary \ref{cor1} that $g \in
C^*_r\left(R\left(\varphi^{l+n}\right)\right)$. Since $f = f\left(v^*\right)^lv^l$
and
$f\left(v^*\right)^l = g$,
this shows that $f$ is in the $C^*$-algebra generated by $v$ and
  $C^*_r\left(R_{\varphi}\right)$. When $f \in C_c\left(\Gamma_{\varphi}(-l,n)\right)$ a similar calculation shows that
$v^lf \in C^*_r\left(R_{\varphi}\right)$ and hence $ f =
\left(v^*\right)^lv^l f$ is in the $C^*$-algebra generated by $v$ and
  $C^*_r\left(R_{\varphi}\right)$.

It follows now from the universal property of
$C^*_r\left(R_{\varphi}\right) \times_{\widehat{\varphi}} \mathbb N$ that there
is a surjective $*$-homomorphism $C^*_r\left(R_{\varphi}\right)
\times_{\widehat{\varphi}} \mathbb N \to C^*_r\left(\Gamma_{\varphi}\right)$
which is the identity on $C^*_r\left(R_{\varphi}\right)$. To show that
it is an isomorphism it suffices, by Proposition 2.1 of \cite{BKR}, to
show that
$$
\left\| \sum_{i \in F} \left(v^*\right)^i a_{i,i} v^i \right\| \leq
\left\| \sum_{i,j \in F} \left(v^*\right)^i a_{i,j} v^j \right\|  
$$
in $C^*_r\left(\Gamma_{\varphi}\right)$ when $F \subseteq \mathbb N$
is a finite set and $\left\{a_{i,j}\right\}_{i,j \in F}$ is any
collection of elements from $C^*_r\left(R_{\varphi}\right)$. This
inequality follows from the existence of the gauge action $\beta$ of $\mathbb
T$ on $C^*_r\left(\Gamma_{\varphi}\right)$ since
$$
 \sum_{i \in F} \left(v^*\right)^i a_{i,i} v^i = \int_{\mathbb T}
 \beta_z \left(  \sum_{i,j \in F} \left(v^*\right)^i a_{i,j} v^j
 \right) dz .
$$
\end{proof}
\end{thm}

\begin{lemma}\label{fullcorner} The endomorphism $\widehat{\varphi} :
  C^*_r\left(R_{\varphi}\right) \to C^*_r\left(R_{\varphi}\right)$ is
  a full corner endomorphism in the sense that the projection
  $\widehat{\varphi}(1) = vv^*$ is full in
  $C^*_r\left(R_{\varphi}\right)$.
\begin{proof} Note that
\begin{equation}\label{eq113}
\widehat{\varphi}(1) (x,y) = vv^*(x,y) = m(x)^{-\frac{1}{2}}
m(y)^{-\frac{1}{2}} 1_{R(\varphi)}(x,y) .
\end{equation}
It follows from Lemma \ref{relation} that the ideal in
$C^*_r\left(R_{\varphi}\right)$ generated by
$\widehat{\varphi}(1)$ contains
$P_{R_{\varphi}}\left(\widehat{\varphi}(1)\right)$. And it follows from
(\ref{eq113}) that $P_{R_{\varphi}}\left(\widehat{\varphi}(1)\right)$ is the
invertible element $m^{-1}$. Hence the ideal in
$C^*_r\left(R_{\varphi}\right)$ generated by
$\widehat{\varphi}(1)$ is all of $C^*_r\left(R_{\varphi}\right)$.
\end{proof}
\end{lemma}

Note that $v$ is a unitary, i.e. $vv^* =1$, if and only if $m =1$
if and only if $\varphi$ is a homeomorphism. In this case
$C^*_r\left(R_{\varphi}\right) = C(X)$ and
$C^*_r\left(\Gamma_{\varphi}\right) \simeq C(X) \times_{\varphi}
\mathbb Z$.
Such crossed products have been intensively studied and we shall have
nothing to add in this case. We therefore restrict attention to the
case where $\varphi$ is surjective, but not injective.

Assume that $\varphi$ is surjective and not injective. Let
$B_{\varphi}$ be the inductive limit of the sequence
\begin{equation}\label{seq}
\begin{xymatrix}{
C^*_r\left(R_{\varphi}\right) \ar[r]^-{\widehat{\varphi}} &
C^*_r\left(R_{\varphi}\right) \ar[r]^-{\widehat{\varphi}} &
C^*_r\left(R_{\varphi}\right) \ar[r]^-{\widehat{\varphi}} &  \dots}
\end{xymatrix}
\end{equation} 
We can then define an automorphism $\widehat{\varphi}_{\infty} \in \Aut
B_{\varphi}$ such that
$\widehat{\varphi}_{\infty} \circ \rho_{\infty,n} = \rho_{\infty,n}
\circ \widehat{\varphi}$,
where $\rho_{\infty,n} : C^*_r\left(R_{\varphi}\right) \to B_{\varphi}$ is the canonical
$*$-homomorphism from the $n$'th level in the sequence (\ref{seq})
into the inductive limit algebra. In this notation the inverse of
$\widehat{\varphi}_{\infty}$ is defined such that
$\widehat{\varphi}_{\infty}^{-1} \circ \rho_{\infty,n} =
\rho_{\infty,n+1}$.

\begin{thm}\label{crossed} Assume that $\varphi$ is surjective and not
  injective. It follows that $p = \rho_{\infty,1}(1) \in B_{\varphi}
  \subseteq B_{\varphi} \times_{\widehat{\varphi}_{\infty}} \mathbb Z$
  is a full projection and
  $C^*_r\left(\Gamma_{\varphi}
\right)$ is $*$-isomorphic to the corner $p\left(B_{\varphi}
  \times_{\widehat{\varphi}_{\infty}} \mathbb Z\right)p$.
\begin{proof} In view of Proposition 3.3. of \cite{St}, which was
  restated in \cite{BKR}, it suffices to check that $\rho_{\infty,1}(1)$
  is a full projection in $B_{\varphi}$. To
  this end it suffices to show that $v^n\left(v^*\right)^n =
  \widehat{\varphi}^n(1)$ a full projection in $C^*_r\left(R_{\varphi}\right)$. By noting that
\begin{equation*}
\begin{split}
&v^n\left(v^*\right)^n(x,y) = \\
&\left[m(x)m(\varphi(x)) \dots
m(\varphi^{n-1}(x))m(y)m(\varphi(y))
\dots m(\varphi^{n-1}(y))\right]^{-\frac{1}{2}}1_{R\left(\varphi^n\right)}(x,y), 
\end{split}
\end{equation*}
this follows from Lemma \ref{relation} as in the proof of Lemma \ref{fullcorner}.
\end{proof}
\end{thm}

\begin{cor}\label{AHcrossed} Assume that $\varphi$ is surjective and
  not injective. It follows that there is a (non-unital) AH-algebra
  $A$ and an automorphism $\alpha$ of $A$ such that
  $C^*_r\left(\Gamma_{\varphi}\right)$ is stably isomorphic to $A
  \times_{\alpha} \mathbb Z$.
\begin{proof} Note that $B_{\varphi}$ is AH
  since $C^*_r\left(R_{\varphi}\right)$ is. Hence the assertion follows
  from Theorem \ref{crossed} by the use of Brown's theorem, \cite{Br}. 
\end{proof}
\end{cor}

One virtue of Theorem \ref{generation} and Theorem \ref{crossed} is
that the crossed product description and the Pimsner-Voiculescu exact
sequence give us a six-term exact
sequence which makes it possible to calculate the $K$-theory of
$C^*_r\left(\Gamma_{\varphi}\right)$ from the $K$-theory of
$C^*_r\left(R_{\varphi}\right)$ and the action of $\widehat{\varphi}$
on $K$-theory. See e.g. \cite{De} and \cite{DS} for such $K$-theory
calculations in the \'etale case.

\subsection{Simplicity of
  $C^*_r\left(\Gamma_{\varphi}\right)$ and $C^*_r\left(R_{\varphi}\right)$}

\begin{thm}\label{simplicity2} Assume that $\varphi$ is
  surjective and not injective. Then $C^*_r\left(\Gamma_{\varphi}\right)$ is simple if
  and only if there is no non-trivial ideal $I$ in
  $C^*_r\left(R_{\varphi}\right)$ such that $\widehat{\varphi} (I)
  \subseteq I$.
\begin{proof} Assume first that $I$ is a non-trivial ideal in
  $C^*_r\left(R_{\varphi}\right)$ which is
  $\widehat{\varphi}$-invariant in the specified way. In the notation
  established before Theorem \ref{crossed}, set $J =
  \overline{\bigcup_n \rho_{\infty,n}(I)}$. Then $J$ is a non-zero
  $\widehat{\varphi}_{\infty}$-invariant ideal in $B_{\varphi}$. To
  prove that $J \neq B_{\varphi}$, we show that $\rho_{\infty,1}(1)
  \notin J$. Indeed, if $\rho_{\infty,1}(1) \in J$ there is, for any
  $\epsilon > 0$, a $k \in \mathbb N$ and an element $a \in I$ such that
 $$
\left\|\widehat{\varphi}^k(1) - a \right\| = \left\|\rho_{\infty,1}(1) -
  \rho_{\infty,k}(a) \right\| \leq \epsilon .
$$
Since $\widehat{\varphi}^k(1)$ is a projection and $I$ is an ideal in
$C^*_r\left(R_{\varphi}\right)$ this implies, with $\epsilon$
appropriately small, that $\widehat{\varphi}^k(1) \in I$. And as
argued in the proof of Theorem \ref{crossed} $\widehat{\varphi}^k(1) =
v^k{v^*}^k$ is a full projection in $C^*_r\left(R_{\varphi}\right)$
and hence $\widehat{\varphi}^k(1) \in I$ implies that $I =
C^*_r\left(R_{\varphi}\right)$, contrary to assumption. Hence $J$ is a
non-trivial ideal in $B_{\varphi}$. Being $\widehat{\varphi}_{\infty}$-invariant it gives rise to a non-trivial ideal in $B_{\varphi}
\times_{\widehat{\varphi}_{\infty}} \mathbb Z$. Since
$C^*_r\left(\Gamma_{\varphi}\right)$ is stably isomorphic to $B_{\varphi}
\times_{\widehat{\varphi}_{\infty}} \mathbb Z$ by Theorem
\ref{crossed} and \cite{Br}, this means that
$C^*_r\left(\Gamma_{\varphi}\right)$ is not simple.

The reversed implication, that $C^*_r\left(\Gamma_{\varphi}\right)$ is
simple when there are no non-trivial $\widehat{\varphi}$-invariant
ideals in $C^*_r\left(R_{\varphi}\right)$ follows from \cite{BKR} (and
\cite{ALNR}), in particular, from Corollary 2.7 of \cite{BKR},
because $C^*_r\left(R_{\varphi}\right)$ is AH, and hence also strongly amenable.

\end{proof}
\end{thm}

\begin{remark}
The inclusion 
$C^*_r\left(R_{\varphi}\right) \subseteq
C^*_r\left(\Gamma_{\varphi}\right)^{\mathbb T}$
is obvious. When $\varphi$ is open the two algebras are identical, but
there are examples of shift spaces where this
is a strict inclusion. Similarly, the inclusion $D_{R_{\varphi}}
\subseteq D_{\Gamma_{\varphi}}$ is obvious and is an identity when
$\varphi$ is open, but a strict inclusion for certain shift spaces. 
\end{remark}

Let $P_{R_{\varphi}} : C^*_r\left(R_{\varphi}\right) \to
D_{R_{\varphi}}$ be the conditional expectation.
Note that
\begin{equation}\label{Pvarphi}
mP_{R_{\varphi}} \circ \widehat{\varphi}(f) =  f \circ \varphi
\end{equation}
for $f \in D_{R_{\varphi}}$. Thus $m P_{R_{\varphi}} \circ
\widehat{\varphi}$ is a unital injective endomorphism
of $D_{R_{\varphi}}$ which we denote by $\overline{\varphi}$.

\begin{lemma}\label{eqlemma} 
\begin{enumerate}
\item[a)] Let $I \subseteq C^*_r\left(R_{\varphi}\right)$ be a
  non-trivial ideal
  in $C^*_r\left(R_{\varphi}\right)$ such that $\widehat{\varphi}(I)
  \subseteq I$. It follows that $J = I \cap D_{R_{\varphi}}$ is
  a non-trivial $R_{\varphi}$-invariant ideal in $D_{R_{\varphi}}$ such that $\overline{\varphi}(J)
  \subseteq J$.
\item[b)] Let $J \subseteq D_{R_{\varphi}}$ be a non-trivial
  $R_{\varphi}$-invariant ideal such that
  $\overline{\varphi}(J) \subseteq J$. It follows that
$$
\widehat{J} = \left\{ a \in C^*_r\left(R_{\varphi}\right) : \ P_{R_{\varphi}}(a^*a)
  \in J \right\}
$$
is a non-trivial ideal in $C^*_r\left(R_{\varphi}\right)$ such that
$\widehat{\varphi}(\widehat{J}) \subseteq \widehat{J}$.
\end{enumerate}
\begin{proof} a): It follows from Lemma \ref{relation} that $P_{R_{\varphi}}(I)
  \subseteq I$ and hence that $I \cap D_{R_{\varphi}}$ is
  not zero since $I$ is not and $P_{R_{\varphi}}$ is faithful. It is not all of
  $D_{R_{\varphi}}$ because it does not contain the
  unit. Finally, it follows that $\overline{\varphi}\left(I \cap
    D_{R_{\varphi}}\right) \subseteq  m P_{R_{\varphi}} \widehat{\varphi}(I)
  \subseteq mP_{R_{\varphi}}(I) \subseteq I \cap D_{R_{\varphi}}$ since
  $\widehat{\varphi}(I) \subseteq I, \ P_{R_{\varphi}}(I) \subseteq I \cap
  D_{R_{\varphi}}$ and $m \in D_{R_{\varphi}}$.

b): Recall that $\widehat{J}$ is a non-trivial ideal by Lemma
\ref{invideal}. Since $P_{R_{\varphi}} \circ \widehat{\varphi} \circ P_{R_{\varphi}} = P_{R_{\varphi}} \circ
\widehat{\varphi}$ it follows that
$$
P_{R_{\varphi}}\left(\widehat{\varphi}(a)^*\widehat{\varphi}(a)\right) = P_{R_{\varphi}} \circ
\widehat{\varphi} \circ P_{R_{\varphi}}(a^*a) = m^{-1}
\overline{\varphi}\left(P_{R_{\varphi}}(a^*a)\right) \in J
$$
when $a \in \widehat{J}$. This shows that
$\widehat{\varphi}\left(\widehat{J}\right) \subseteq \widehat{J}$. 

 \end{proof}
\end{lemma}

\begin{thm}\label{simplicityD} Assume that $\varphi$ is
  surjective and not injective. Then $C^*_r\left(\Gamma_{\varphi}\right)$ is simple if
  and only if there is no non-trivial $R_{\varphi}$-invariant ideal $J$ in
  $D_{R_{\varphi}}$ such that $\overline{\varphi} (J)
  \subseteq J$.
\begin{proof} Combine Lemma \ref{eqlemma} and Theorem \ref{simplicity2}.
\end{proof}
\end{thm}

When $\varphi$ is open, and hence a
  local homeomorphism, Theorem \ref{simplicityD} follows from
  Proposition 4.3 of \cite{DS}.

For all $k,l \in \mathbb N$, set
$$
X_{k,l}= \left\{ x \in X : \ \# \varphi^{-k}\left(\varphi^k(x)\right) =
  l \right\} .
$$

\begin{thm}\label{ultsimplthm} Assume that $\varphi$ is surjective and
  not injective. The following conditions are equivalent:
\begin{enumerate}
\item[1)] $C^*_r\left(\Gamma_{\varphi}\right)$ is simple.
\item[2)] For every open subset $U \subseteq X$
  and $k, l \in
\mathbb N$ such that $U \cap X_{k,l} \neq \emptyset$
there is an $m \in \mathbb N$ such that
$$
\bigcup_{j=0}^{m} \varphi^{j}\left(U \cap X_{k,l}\right) = X.
$$
\end{enumerate}
\begin{proof} For the organization of the proof it is convenient to
  observe that condition 2) is equivalent to the following:
\begin{enumerate}
\item[2')] For every open subset $U \subseteq X$
  and $k, l \in
\mathbb N$ such that $U \cap X_{k,l} \neq \emptyset$
there is a $m \in \mathbb N$ such that
$$
\varphi^{m+k}\left(\bigcup_{j=0}^{m-1} \varphi^{-j}\left(U \cap X_{k,l}\right)\right) = X.
$$
\end{enumerate}

1) $\Rightarrow$ 2'): Assume first that $C^*_r\left(\Gamma_{\varphi}\right)$ is
  simple. If condition 2') fails there is an open subset $U$
  in $X$ and a pair $k,l \in \mathbb N$ such that $U \cap X_{k,l} \neq
  \emptyset$ and for each $m \in \mathbb N$ there is an element $x_m \in X$ such
  that 
\begin{equation}\label{A17}
\varphi^{-m-k}\left(x_m\right) \cap
\left(\bigcup_{j=0}^{m-1} \varphi^{-j}\left(U \cap X_{k,l} \right)\right) =
  \emptyset .
\end{equation}
Let $R\left(\varphi^k\right)|_U$ be the reduction of
$R\left(\varphi^k\right)$ to $U$, i.e. 
$$R\left(\varphi^k\right)|_U =
\left\{ (x,y) \in R\left(\varphi^k\right) : \ x,y \in U\right\},
$$
which is an open subgroupoid of $R\left(\varphi^k\right)$. It follows
from Lemma \ref{openincl} that there is an isometric inclusion
$C^*_r\left(R\left(\varphi^k\right)|_U \right) \subseteq
C^*_r\left(R\left(\varphi^k\right)\right)$. In the notation used in
Lemma \ref{iotadense}, let $\psi_{x_m}$ be the irreducible
representation of $C^*_r\left(R\left(\varphi^{k+m}\right)\right)$ corresponding to
$x_m$. Let $\{p_j\}$ be the central
projections of Lemma \ref{centralproj} relative to $\varphi^k$. It follows from (\ref{A17}) that $\psi_{x_m}\left(
   \widehat{\varphi}^{j}
   \left(C^*_r\left(R\left(\varphi^{k}\right)|_U\right)p_l\right)\right) = 0$ for all $j \in \left\{0,1,2, \dots, m-1\right\}$. By composing the normalized trace of $M_{\varphi^{-m-k}(x_m)}(\mathbb C)$ with $\psi_{x_m}$ we obtain in this way, for each $m \in \mathbb N$, a trace state $\omega_m$ on $C^*_r\left(R\left(\varphi^{m+k}\right)\right)$
which annihilates
$\widehat{\varphi}^j\left(C^*_r\left(R\left(\varphi^k\right)|_U\right)p_l\right)$
for $j = 0,1,\dots,m-1$. For each $m$ we choose a state extension
$\omega'_m$ of $\omega_m$ to $C^*_r\left(R_{\varphi}\right)$. Any weak*
condensation point of the sequence $\{\omega'_m\}$ will be a trace state $\omega$ on
$C^*_r\left(R_{\varphi}\right)$ which annihilates
$\widehat{\varphi}^j\left(C^*_r\left(R\left(\varphi^k\right)|_U\right)p_l\right)$
for all $j \in \mathbb N$. It follows that the closed two-sided ideal $I$
in $C^*_r\left(R_{\varphi}\right)$ generated by
$$
\bigcup_{j=0}^{\infty} \widehat{\varphi}^j\left(C^*_r\left(R\left(\varphi^k\right)|_U\right)p_l\right)
$$
is contained in $\left\{a \in C^*_r\left(R_{\varphi}\right) : \
  \omega(a^*a) = 0\right\}$. Hence $I$ is a non-trivial ideal in
$C^*_r\left(R_{\varphi}\right)$. Since $\widehat{\varphi}(I) \subseteq
I$ this
contradicts the simplicity of $C^*_r\left(\Gamma_{\varphi}\right)$ by
Theorem \ref{simplicity2}.

 2') $\Rightarrow$ 1): Let $I \subseteq
C^*_r\left(R_{\varphi}\right)$ be a non-zero closed twosided ideal such that
$\widehat{\varphi}(I) \subseteq I$. By Theorem \ref{simplicity2} it
suffices to show that $1 \in I$. Since $I \neq
0$ there is a $k \in \mathbb N$ such that $I \cap C^*_r\left(R\left(\varphi^k\right)\right) \neq
0$. Let $d' \in I \cap C^*_r\left(R\left(\varphi^k\right)\right)$ be an element
with $\|d'\| = 1$. There
is then an $l$ such that $d'p_l \neq 0$. In particular, there is an $a \in X$ such that
$\psi_a\left(d'p_l\right) \neq 0$.
Since $I \cap C^*_r\left(R\left(\varphi^k\right)\right)$ is an ideal
in $C^*_r\left(R\left(\varphi^k\right)\right)$, $p_l$ is central in
$C^*_r\left(R\left(\varphi^k\right)\right)$ and $\psi_a
\left(C^*_r\left(R\left(\varphi^k\right)\right)\right) \simeq
M_{\varphi^{-k}(a)}(\mathbb C)$, there is a positive element $d \in I \cap
C^*_r\left(R\left(\varphi^k\right)\right)$ such that $dp_l(\xi , \xi) >1$
for all $\xi \in \varphi^{-k}(a) \cap X_{k,l}$. It follows from Lemma \ref{contonlj} that the map
$\xi  \mapsto dp_l(\xi, \xi)$
is continuous on $X_{k,l}$. 
There is therefore an open set
$W$ in $X$ such that $W\cap X_{k,l} \supseteq \varphi^{-k}(a) \cap
X_{k,l}$ and
\begin{equation}\label{simpleq1}
dp_l(\xi ,\xi)  \geq 1 .
\end{equation}
for all $\xi \in W \cap X_{k,l}$. Since condition 2') holds there
is an $m \in \mathbb N$ such that
\begin{equation}\label{cruxAA}
\varphi^{m+k}\left(\bigcup_{j=0}^{m-1} \varphi^{-j}\left(W \cap X_{k,l}\right)\right) = X .
\end{equation}
 Now, assuming that $1 \notin I$, there
is an irreducible ideal in $C^*_r\left(R\left(\varphi^{m+k}\right)\right)$
  which contains $I \cap C^*_r\left(R\left(\varphi^{m+k}\right)\right)$. It
  follows then from Lemma \ref{iotadense} and the $\widehat{\varphi}$-invariance
  of $I$ that there is a point $b \in
  X$ such that 
\begin{equation}\label{simpleq2}
\left\|\psi_b\left(\widehat{\varphi}^i\left(dp_l\right)\right)\right\| <
  \inf_{y \in X, j \in \left\{0,1,2, \dots, m-1\right\}} \
\left(m(y)m\left(\varphi(y)\right) \cdots m\left(\varphi^{j-1}(y)\right)\right)^{-1}
\end{equation}  
for all $i \in \{0,1,2, \dots, m-1\}$. It follows from (\ref{cruxAA})
that there is a $j_0 \in
\{0,1,2,\dots, m-1\}$ and an element $z \in \varphi^{-m-k}(b) \cap \varphi^{-j_0}\left(W \cap X_{k,l}\right)$. Then
\begin{equation}\label{jameslast}
\begin{split}
& \left\|\psi_b \circ \widehat{\varphi}^{j_0}\left(dp_l\right)\right\|
 \geq \widehat{\varphi}^{j_0}\left(dp_l\right)(z,z) \\
& =
 \left(m(z)m\left(\varphi(z)\right) m\left(\varphi^2(z)\right) \dots
   m\left(\varphi^{j_0-1}(z)\right)\right)^{-1}dp_l\left(\varphi^{j_0}(z),\varphi^{j_0}(z)\right) \\
& \geq \left(m(z)m\left(\varphi(z)\right) m\left(\varphi^2(z)\right) \dots
   m\left(\varphi^{j_0-1}(z)\right)\right)^{-1}, 
\end{split}
\end{equation}
where we in the last step used that $\varphi^{j_0}(z) \in W \cap
X_{k,l}$ so that (\ref{simpleq1}) applies. (\ref{jameslast}) contradicts (\ref{simpleq2}).
\end{proof}
\end{thm}

It is easy to modify the proof of Theorem \ref{ultsimplthm} to obtain
the following

\begin{thm}\label{relsimplealg} Assume that $\varphi$ is surjective and
  not injective. The following conditions are equivalent:
\begin{enumerate}
\item[1)] $C^*_r\left(R_{\varphi}\right)$ is simple.
\item[2)] For every open subset $U \subseteq X$
  and $k, l \in
\mathbb N$ such that $U \cap X_{k,l} \neq \emptyset$
there is a $m \in \mathbb N$ such that
$$
 \varphi^{m}\left(U \cap X_{k,l}\right) = X.
$$
\end{enumerate}

\end{thm}

By Theorem \ref{simplicity1} the two conditions, 1) and 2), in Theorem
\ref{relsimplealg} are equivalent to the absence of any non-trivial
$R_{\varphi}$-invariant ideal in $D_{R_{\varphi}}$.

When $\varphi$ is open and hence a local homeomorphism it was observed
in \cite{A} that the sets
$X_{k,l}$ are all open by a result of Eilenberg. So in this case
condition 2) of Theorem \ref{ultsimplthm} is equivalent to strong
transitivity of $\varphi$ in the sense of \cite{DS}: For every
non-empty open set $U$ of $X$ there is an $m \in \mathbb N$ such that
$\bigcup_{j=0}^m \varphi^j(U) = X$. Similarly when $\varphi$ is open condition 2) of Theorem
\ref{relsimplealg} is equivalent to exactness of $\varphi$: For every
non-empty open set $U$ of $X$ there is an $m \in \mathbb N$ such that
$\varphi^m(U) = X$.  
So when $\varphi$ is a local
homeomorphism Theorem \ref{ultsimplthm} follows from Proposition 4.3
of \cite{DS} and Theorem \ref{relsimplealg} from Proposition 4.1 of
\cite{DS}.

We show next that it is possible to use the methods of
this paper to improve the known simplicity criteria in the \'etale
case to handle a non-surjective local homeomorphism of a locally
compact space. Let $X$ be a locally compact Hausdorff space and $\varphi : X \to X$
a local homeomorphism.
We say that $\varphi$ is \emph{irreducible} when
\begin{equation}\label{transicond}
X = \bigcup_{0 \leq i, j } \varphi^{-i}\left(\varphi^j(U)\right) .
\end{equation}
for every open non-empty set $U$ in $X$. As observed in \cite{EV} a simple
argument shows that $\varphi$
  is irreducible if and only if there is no non-trivial open subset $V
  \subseteq X$ such that $\varphi^{-1}(V) =V$.

\begin{thm}\label{localhomeomorphism} Let $X$ be a locally compact
  second countable
  Hausdorff space and $\varphi : X \to X$ a local
  homeomorphism. The following are
  equivalent:
\begin{enumerate}
\item[1)] $C^*_r(\Gamma_{\varphi})$ is simple.
\item[2)] $\left\{x \in X : \
  \varphi^k(x) = x\right\}$ has empty interior for each $k \geq 1$ and
$\varphi$ is irreducible.
\item[3)] There is a point in $X$ which is not pre-periodic under $\varphi$ and
$\varphi$ is irreducible. 
\end{enumerate}
\begin{proof} 1) $\Rightarrow$ 2): Assume that $\left\{x \in X : \
  \varphi^k(x) = x\right\}$ contains a non-empty open set $V$ for some $k \geq
1$. Set $W = \bigcup_{j=0}^{k-1}\varphi^j(V)$. Then the reduction 
$$
\Gamma_{\varphi}|_W = \left\{(x,k,y) \in \Gamma_{\varphi} : \ x,y \in W \right\}
$$
is an \'etale groupoid in itself and $C^*_r\left(\Gamma_{\varphi}|_W\right)$ is
a $C^*$-subalgebra of $C^*_r(\Gamma_{\varphi})$. It is easy to check that
$$
C^*_r\left(\Gamma_{\varphi}|_W\right) = \overline{C_0(W)C^*_r(\Gamma_{\varphi})C_0(W)},
$$
showing that $C^*_r\left(\Gamma_{\varphi}|_W\right)$ is a hereditary
$C^*$-subalgebra. Note that $C^*_r(\Gamma_{\varphi})$ is separable since we assume that
$X$ is second countable. Since we assume that
$C^*_r\left(\Gamma_{\varphi}\right)$ is simple we can then apply \cite{Br} to conclude
that $C^*_r(\Gamma_{\varphi})$ is stably isomorphic to
$C^*_r\left(\Gamma_{\varphi}|_W\right)$. However, since $\varphi$ is $k$-periodic
on $W$, every orbit of an element in $W$ is a $\Gamma_{\varphi}|_W$-invariant
closed subset of $W$. As $C^*_r\left(\Gamma_{\varphi}|_W\right)$ must be
simple since $C^*_r(\Gamma_{\varphi})$ is, it follows from Corollary
\ref{idealcor1} and (the proof of) Corollary \ref{invopne} that $W$
must be a single orbit. But then
$$
C^*_r\left(\Gamma_{\varphi}|_W\right) \ \simeq \ C(\mathbb T)
\otimes M_{k'}(\mathbb C) 
$$
where $k' \leq k$ is the number of elements in $W$. This algebra is obviously not simple, contradicting the assumption
that $C^*_r(\Gamma_{\varphi})$ is. It follows that $\left\{x \in X : \
  \varphi^k(x) = x\right\}$ must have empty interior for each $k \geq
1$. 

It follows from Corollary \ref{idealcor1} and (the proof of) Corollary
\ref{invopne} that $X$ contains no non-trivial open $\Gamma_{\varphi}$-invariant
subset. It is easy to see that this is equivalent to the assertion
that (\ref{transicond}) holds for every non-empty open subset $U$.

2) $\Rightarrow$ 3): Assume to reach a contradiction that every element
of $X$ is pre-periodic under $\varphi$. This
means that
\begin{equation}\label{unioneq}
X \  = \bigcup_{n \geq 1, \ k \geq 0} \varphi^{-k}\left(\Per_n X\right)
\end{equation}
where $\Per_n X = \left\{y \in X : \ \varphi^n(y) = y\right\}$. It
follows from the Baire category theorem that there are $n \geq 1$, $k
\geq 0$ and a non-empty open set $V \subseteq \varphi^{-k}\left(\Per_n
  X\right)$. Since $\varphi$ is open this implies that $\varphi^k(V)$ is
an open subset of $\Per_n X$, contradicting our assumption.

3) $\Rightarrow$ 1): As we observed above irreducibility of
$\varphi$ is equivalent to the absence of any non-trivial
$\Gamma_{\varphi}$-invariant open subset in $X$. Furthermore, a point $x$ of $X$
which is not pre-periodic under $\varphi$ must have trivial isotropy
group in $\Gamma_{\varphi}$. Hence the simplicity of
$C^*_r(\Gamma_{\varphi})$ follows from Theorem \ref{simplicity1}.
\end{proof}
\end{thm}

Concerning the simplicity of $C^*_r\left(R_{\varphi}\right)$ when
$\varphi$ is open we get the following
conclusion directly from Corollary \ref{invopne}. It generalizes
Proposition 4.1 of \cite{DS}.

\begin{thm}\label{reletalesimpl} Let $X$ be a locally compact
  Hausdorff space and $\varphi: X \to X$ a local homeomorphism. It
  follows that $C^*_r\left(R_{\varphi}\right)$ is simple if and only
  if $\bigcup_{k=0}^{\infty} \varphi^{-k}\left(\varphi^k(U)\right) =
  X$ for every open non-empty subset $U \subseteq X$.
\end{thm}

\subsection{Subshifts: Carlsen-Matsumoto algebras}\label{sec3} K. Matsumoto was the first to encode structures from general subshifts
in a $C^*$-algebra, \cite{Ma1}, generalizing the original construction of Cuntz
and Krieger, \cite{CK}. Later, slightly different constructions were
suggested by Carlsen and Matsumoto, \cite{CM}, and by Carlsen,
\cite{C}. The exact relation between the various constructions is a
little obscure. Some of the known connections between them are
described in \cite{CM} and \cite{CS}. As we shall see the
approach we take here, based on the groupoids of Renault,
Deaconu and Anantharaman-Delaroche, gives rise to the algebras introduced by
Carlsen in \cite{C}.

Set $A = \{1,2, \dots, n\}$ and $A^{\mathbb N} =
\left\{\left(x_1,x_2,x_3, \dots\right) : \  x_i \in A
\right\}$. We consider $A^{\mathbb N}$ as a compact metric space with the
metric
$$
d(x,y) = \sum_{i=1}^{\infty} 2^{-i} \left|x_i - y_i\right| .
$$ 
The shift $\sigma$ acts on $A^{\mathbb N}$ in the usual way:
$\sigma(x)_i = x_{i+1}$. Let $S \subseteq A^{\mathbb N}$ be a subshift, i.e. $S$
is closed and 
$\sigma(S) = S$.  
Such a subshift defines in a canonical way an
abstract language whose words $\mathbb W(S)$ are the finite strings of
'letters' from the 'alphabet' $A$ which occur in an element from
$S$. We refer to \cite{LM} for more on subshifts. 

Since $\sigma : S \to S$ is locally injective we can apply the
construction of the previous section to obtain a semi \'etale groupoid
which we denote by $\Gamma_S$. Similarly, the corresponding semi
\'etale equivalence relation will be denoted by $R_S$ in this
setting. Given a word $u  \in \mathbb W(S)$ of length $|u| =n$, set
$$
C(u) = \left\{ x \in S :  \ x_1x_2\dots x_n =u \right\} .
$$
These are the \emph{standard cylinder} sets in $S$ and they form a base for
the topology. Now set
$$
C(u,v) =  C(u) \cap \sigma^{-|u|}\left(\sigma^{|v|}\left(C(v)\right)\right) .
$$
Thus $C(u,v)$ consists of the elements of $C(u)$ with the property
that when the prefix $u$ is replaced by $v$ the infinite row of letters is
still an element of $S$. Since the empty word $\emptyset$ by
convention is also a word in $\mathbb W(S)$, with cylinder
$C_{\emptyset} = S$, we have that $C(u) = C(u,\emptyset)$. 
While the cylinder sets are both closed and
open, the set $C(u,v)$ is in general only closed. The characteristic
functions $1_{C(u,v)},  u,v \in \mathbb W(S)$, generate a separable
$C^*$-subalgebra in $l^{\infty}(S)$ which we denote by $D_S$. The
$C^*$-algebra $\mathcal O_S$ of Carlsen from \cite{C} is generated by partial
isometries $s_u, u \in \mathbb W(S)$, such that 
$$
s_us_v  = s_{uv}
$$
when $uv \in \mathbb W(S)$, $s_us_v = 0$ when $uv \notin \mathbb W(S)$, and such that 
$s_vs_u^*s_us_v^*$,
$u,v \in \mathbb
W(S)$, are projections which generate an abelian $C^*$-subalgebra isomorphic to
$D_S$ under a map sending $s_vs_u^*s_us_v^*$ to $1_{C(v,u)}$,
cf. \cite{CS}.

The algebra $\mathcal O_S$ is blessed with the
following universal property
which enhances its applicability: 

\smallskip

\emph{(A) When $B$ is a $C^*$-algebra
containing partial isometries $S_u, u \in \mathbb W(S)$, such that
$$
S_{uv} = S_uS_v
$$ 
when $uv \in \mathbb W(S)$, $S_uS_v = 0$ when $uv \notin \mathbb W(S)$, and admitting a
$*$-homomorphism $D_S \to B$ sending $1_{C(v,u)}$ to
$S_vS_u^*S_uS_v^*$ for all $u,v \in \mathbb W(S)$, then there is a $*$-homomorphism $\mathcal O_S \to
B$ sending $s_u$ to $S_u$ for all $u \in \mathbb W(S)$.}

\smallskip

In particular, it follows from (A) that there is a continuous action $\gamma$ of
the unit circle $\mathbb T$ on $\mathcal O_S$ such that $\gamma_z(s_u)
= z^{|u|}s_u$ for all $z \in \mathbb T$ and all $u \in \mathbb W(S)$. This action is
called the \emph{gauge action}, cf. \cite{C},\cite{CS}.

The universal property (A) of $\mathcal O_S$ is established in Theorem 10 of
\cite{CS}. As we shall show in the following proof, property (A) will provide us
with a $*$-homomorphism $\mathcal O_S \to C^*_r\left(\Gamma _S\right)$
which we show is surjective. To conclude that it is also injective we
shall need that $\mathcal O_S$ is a crossed product of a type dealt with
by Exel and Vershik in \cite{EV}. In the notation of \cite{CS} and
\cite{EV}, $\mathcal O_S = D_S \rtimes_{\alpha, \mathcal L} \mathbb
N$, where $\alpha$ is the endomorphism of $D_S$ defined such that
$\alpha(f) = f \circ \sigma$,
and the transfer operator $\mathcal L$ is given by
$$
\mathcal L(f)(x) = \frac{1}{\# \sigma^{-1}(x)} \sum_{y \in
  \sigma^{-1}(x)} f(y),
$$
cf. \cite{CS}. Note that both $\sigma$ and $\mathcal L$ are
unital and faithful so that the Hypotheses 3.1
of \cite{EV} are satisfied. Furthermore, by inspection of the proof of
Theorem 10 in \cite{CS} one sees that the gauge action $\gamma$
of $\mathcal O_S$ is the same as the gauge action considered in
\cite{EV}. In this way we can use Theorem 4.2 of
\cite{EV} to supplement property (A) with the following 'gauge
invariant uniqueness property':

\smallskip

\emph{(B) Let $B$ be $C^*$-algebra and $\lambda : \mathcal O_S \to B$ a
$*$-homomorphism which is injective on $D_S$. Assume that $B$ admits
a continuous action of $\mathbb T$ by automorphisms such that
$\lambda$ is equivariant with respect to the gauge-action on $\mathcal
O_S$. It follows that $\lambda$ is injective.}

 \smallskip

\begin{thm}\label{matsumoto} Let $S$ be a one-sided subshift. Then the
  Carlsen-Matsumoto algebra $\mathcal O_S$ 
   is $*$-isomorphic to the $C^*$-algebra
  $C^*_r\left(\Gamma_S\right)$ under a $*$-isomorphism which maps
  $D_S$ onto $D_{\Gamma_S}$.
 \begin{proof} When $u \in \mathbb W(S)$ is a word, we let $t_u \in
  B_c(\Gamma_S)$ be the
  characteristic function of the set
$$
\left\{ (x, l, y) \in S \times \mathbb Z \times S : \ x \in C(u), l =
  |u|, \ y_i = x_{|u|+i}, \ i \geq 1\right\}.
$$
Note that $\left\{ (x, l, y) \in S \times \mathbb Z \times S : \ x \in C(u), l =
  |u|, \ y_i = x_{|u|+i}, \ i \geq 1\right\}$ is an open and compact subset of
${\Gamma_S}\left(|u|,0\right)$ and hence of
$\Gamma_S$. Therefore $t_u \in
C_c\left(\Gamma_S\right)$. Straightforward calculations show that 
\begin{equation}\label{oz3}
t_vt_u^*t_ut_v^* = 1_{C(v,u)}
\end{equation}
for all $u,v \in \mathbb W(S)$ when we identify $S$ with the unit
space of $\Gamma_S$, and that
$t_ut_v = t_{uv}$
when $uv \in \mathbb W(S)$, and $t_ut_v = 0$ when $uv \notin \mathbb W(S)$. It follows therefore from the universal
property (A)
of $\mathcal O_S$ described above that there is a
$*$-homomorphism $\lambda : \mathcal O_S \to C^*_r(\Gamma_S)$ such that
$\lambda\left(s_u\right) = t_u$ for all $u \in \mathbb W(S)$.

To see that $\lambda$ is surjective note first that it follows from
(\ref{oz3})  that $1_{C(u)} = 1_{C(u,u)}$ is in the range of $\lambda$ for alle $u
\in \mathbb W(S)$. Note next that $t_ut_v^* = 1_{A(u,v)}$
where 
\begin{equation}\label{standardn}
\begin{split}
&A(u,v) = \\
&\left\{ (x,k,y) \in S \times \mathbb Z \times S : \ k = |u| -
  |v|, \ x \in C(u), \ y \in C(v), \ x_{|u|+i} = y_{|v|+i}, \ i \geq
  1\right\}
\end{split}
\end{equation}
is a compact and open subset of $\Gamma_S$. In fact, sets
of this form constitute a base for the topology of $\Gamma_S$ so in
order to show that every element of $C_c\left(\Gamma_S\right)$ is in
the range of $\lambda$, which implies that $\lambda$ is
surjective, it will be enough to consider an $f
\in C_c\left(\Gamma_S\right)$ such that $f$ has
support in $A(u,v)$ for some $u,v \in \mathbb W(S)$. Let $\epsilon >
0$. Note that the
range map $r$ is
injective on $A(u,v)$. By combining this fact with the continuity of
$f$ it follows that we can find words $u_i, i = 1,2,
\dots, N$, in $\mathbb W(S)$ such that $C(u) = \bigcup_{i=1}^N C(u_i)$ and
\begin{equation}\label{eq11}
\left| f(\xi) - f(\xi')\right| \leq \epsilon
\end{equation}
when
$\xi, \xi' \in A\left(u,v\right) \cap r^{-1}\left(C\left(u_i\right)\right)$, and
such that 
\begin{equation}\label{oz2}C\left(u_i\right) \cap
C\left(u_j\right) = \emptyset, i \neq j.
\end{equation} 
Define a function $h : S \to \mathbb
C$ such that
$$
h(x) = \begin{cases} f(x,|u|-|v|,y), \ & \ \text{when} \ (x,|u|-|v|,y) \in
  A(u,v) \ \text{for some} \ y \in C(v) \\ 0, \ & \ \text{otherwise} .\end{cases}
$$   
Then $h$ is bounded and supported in $C(u)$. Let $J = \left\{ i : \
  A\left(u,v\right) \cap r^{-1}\left(C(u_i)\right)  \neq \emptyset \right\}$. For each
$i \in J$ we pick an element $\xi_i \in A\left(u,v\right) \cap r^{-1}\left(C(u_i)\right)$ and define $k : S \to \mathbb
C$ such that
$$
k = \sum_{i \in J} h\left(\xi_i\right)1_{C(u_i)} .
$$
Then $k,h$ are both bounded and compactly supported in
$S$ and $k$ is in the range of $\lambda$. Furthermore, it follows from
(\ref{eq11}) and (\ref{oz2}) that
\begin{equation}\label{eq12}
\left\| h-k\right\|_{\infty} \leq \epsilon 
\end{equation}
in $B_c\left(\Gamma_S^0\right)$. By using the
canonical inclusion $B_c\left(\Gamma_S^0\right) \subseteq
B_c\left(\Gamma_S\right)$ we consider $h$ and $k$ as elements of
$B_c\left(\Gamma_S\right)$, and find then that 
$f = h
\star t_{u}t_{v}^*$. Hence
\begin{equation}\label{eq13}
\left\|f - k \star t_ut_v^*\right\|_{C^*_r\left(\Gamma_S\right)} \leq \|h-k\|_{\infty} \left\|t_{u}t_{v}^*\right\|_{C^*_r\left(\Gamma_S\right)} .
\end{equation}
It follows from (\ref{oz3}) that $t_ut_v^*$ is a partial isometry and hence
that $\left\|t_ut_v^*\right\| \leq 1$. We can therefore combine
(\ref{eq13}) and (\ref{eq12}) to get the estimate
$$
\left\|f - k \star t_ut_v^*\right\| \leq \epsilon.
$$ 
Since $\epsilon > 0$ is arbitrary and
  $k \star t_ut_v^*$ is in the range of
  $\lambda$, it follows that so is $f$, completing the proof that
  $\lambda$ is surjective. Note that $1_{C\left(u_i\right)} \star
  t_ut_v^* = 1_{A(u_i,v_i')}$ for an appropriate word $v'_i \in
  \mathbb W(S)$ so that $k \star t_ut_v^*$ is a linear combination of
  such characteristic functions.

To see that $\lambda$ is injective note first that (\ref{oz3}) implies
that $\lambda$ is injective on $D_S$. Therefore property (B) above
shows that the injectivity of $\lambda$ will follow
if we can exhibit a continuous action $\beta : \mathbb T \to \Aut
C^*_r\left(\Gamma_S\right)$ such that $\beta_z\left(t_u\right) =
z^{|u|}t_u$ for all $u \in \mathbb W(S)$. The gauge action from Lemma
\ref{gaugeaction} is such an
action.

It remains to show that 
\begin{equation}\label{oz5}
\lambda\left(D_S\right) = D_{\Gamma_S}.
\end{equation} 
The inclusion $\lambda\left(D_S\right) \subseteq D_{\Gamma_S}$ follows
from (\ref{oz3}) and the definition of $\lambda$, so it remains only
to show that $D_{\Gamma_S} \subseteq \lambda\left(D_S\right)$. To
this end we let $1_{A(u,v)}$ denote the characteristic
function of the set (\ref{standardn}). Let $u \in \mathbb W(S)$, and
let $F \subseteq \mathbb W(S)$ be a finite set of words in $S$, not necessarily of the same length as $u$. We set
\begin{equation}\label{skew}
C'(u; F) = \left\{ x \in C(u) : \ \forall v \in F \  \exists y^v \in C(v)
  \ \text{such that} \ y^v_{|v| +i} = x_{|u| +i} , \ i \geq 1 \right\}
.
\end{equation}
Thus $C'(u;F) = \bigcap_{v \in F} C(u,v)$.
Let $v_1,v_2, \dots, v_N$ be the elements of $F$. It is
straightforward to check that
\begin{equation}\label{form200}
1_{C'(u;F)} = P_{\Gamma_S}\left(1_{A(u,v_1)} \star 1_{A(v_1,v_2)} \star \dots
  \star 1_{A\left(v_{N-1},v_N\right)} \star 1_{A(v_N,u)} \right)
\end{equation}
where $P_{\Gamma_S} : C^*_r\left(\Gamma_S\right) \to D_{\Gamma_S}$ is the
conditional expectation of Lemma \ref{condexp} corresponding to $\Gamma_S$.
In particular, $1_{C'(u;F)} \in D_{\Gamma_S}$. As we have just shown
  every element of $C_c(\Gamma_S)$ can be approximated arbitrarily
  well in $C^*_r(\Gamma_S)$ by a linear combination of functions of the form
  $1_{A(u,v)}$. It follows that $C^*_r(\Gamma_S)$ is the closed linear
  span of elements of the form
\begin{equation}\label{weq1}
1_{A(u_1,v_1)} \star 1_{A(u_2,v_2)} \star \dots \star 1_{A(u_N,v_N)} .
\end{equation}
By using that
$1_{A(u,v)} = \sum_{i=1}^n 1_{A(ui,vi)}$
we can write the convolution product (\ref{weq1}) as a sum of similar
products, with the additional property that $|v_i| = |u_{i+1}|,
i=1,2,\dots, N-1$. Then 
$P_{\Gamma_S}\left(1_{A(u_1,v_1)} \star 1_{A(u_2,v_2)}
  \star \dots \star 1_{A(u_N,v_N)}\right) = 0$ 
unless $v_1 = u_2, v_2 = u_3, \dots, v_N = u_1$, in
which case
$P_{\Gamma_S}\left(1_{A(u_1,v_1)} \star 1_{A(u_2,v_2)}
  \star \dots \star 1_{A(u_N,v_N)}\right) = 1_{C'(u_1;F)}$ where $F =
\left\{v_1,v_2, \dots, v_{N-1}\right\}$. This shows that
$D_{\Gamma_S}$ is the closed linear span of projections of the form
$1_{C'(u;F)}$ for some $u$ and $F$. Since 
$1_{C'(u;F)} = 1_{C(u,v_1)}1_{C(u,v_2)}\cdots
1_{C\left(u,v_k\right)}$ when $F = \left\{v_1,v_2, \dots, v_k\right\}$ and since
$\lambda\left(s_us_{v_i}^*s_{v_i}s_u^*\right) = 1_{C\left(u,v_i\right)}$,
we obtain (\ref{oz5}).
\end{proof}
\end{thm}

Let $u \in \mathbb W(S)$, and
let $F$ be a finite set of words in $S$ of the same length as $u$. We set
$$
C(u; F) = \left\{ x \in C(u) : \ \forall v \in F \  \exists y^v \in C(v)
  \ \text{such that} \ y^v_{|v| +i} = x_{|u| +i} , \ i \geq 1 \right\}
.
$$
Similar sets were used in the proof of Theorem \ref{matsumoto}, but
note that we now require the words in $F$ to have the same length as
$u$. We will then call $C(u;F)$ 
\emph{a generalized cylinder} in $S$. Following the notation used in Theorem \ref{ultsimplthm} we set
$$
S_{k,l} = \left\{x \in S: \ \# \sigma^{-k}\left(\sigma^k(x)\right) =
  l\right\}
$$
for all $k,l \in \mathbb N$.

\begin{lemma}\label{sss} a) When $U \subseteq S$ is an open subset and
  $k,l \in \mathbb N$ are numbers such that $U \cap S_{k,l} \neq
  \emptyset$ there is a non-empty generalized cylinder $C(u;F)$ such
  that $C(u;F) \subseteq U \cap S_{k,l}$.

b) When $C(u;F)$ is a non-empty generalized cylinder there is an open
subset $U \subseteq S$ and natural numbers $k,l \in \mathbb N$ such that $U \cap S_{k,l}
\neq \emptyset$ and $U \cap S_{k,l} \subseteq C(u;F)$.
\begin{proof} a) Let $x \in U \cap S_{k,l}$, and let $\mathbb W_k(S)$
  denote the set of words in $S$ of length $k$. There are then exactly
  $l$ words $v_1,v_2,\dots, v_l$ in $\mathbb W_k(S)$ such that $v_jx_{[k+1,\infty)} =
  v_jx_{k+1}x_{k+2}x_{k+3} \dots \ \in S$. The word $x_1x_2\dots x_k$
  is one of them and we arrange that it is $v_1$. For each $j \in
  \{1,2,\dots, l\}$ and some $m > k$, set $
w_j = v_jx_{k+1}x_{k+2} \dots x_m$.
We can then choose $m$ so large that
$x \in C\left(w_1; \left\{w_2,w_3, \dots, w_l\right\}\right) \subseteq
U \cap S_{k,l}$.

b) Let $F'$ be a maximal collection of words from $\mathbb W_{|u|}(S)$
with the property that $F \subseteq F'$ and $C(u;F') \neq
\emptyset$. Let $x\in C(u;F')$ and set $k = |u|, l =  \# F'$. For each
word $v \in \mathbb W_k(S) \backslash F'$ there is a natural number
$m_v$ such that
$$
vx_{\left[k+1,i\right]} \notin \mathbb W(S)
$$
when $i \geq m_v$. Set $m = \max_v m_v$. Then
$$
x \in S_{k,l} \cap C\left(x_{[1,m]}\right) \subseteq C(u;F') \subseteq C(u;F).
$$
Let $U = C\left(x_{[1,m]}\right)$. This handles the case when $F' \neq
\mathbb W_k(S)$. In case $F' = \mathbb W_k(S)$ we have that $l = \# \mathbb W_k(S)$,
and we can then take $U = C\left(x_{[1,k]}\right)$. 
 
\end{proof}
\end{lemma}

\begin{thm}\label{subshiftsimple} Let $S$ be an infinite
  one-sided subshift. Then $C^*_r\left(\Gamma_S\right)$ is simple if
  and only if
  the following holds: For every non-empty generalized cylinder
  $C(u;F)$ there is an $m \in \mathbb N$ with the property that for
  all $x \in S$ there is an element $y \in C(u;F)$ and a $k \in
  \{0,1,2,\dots, m\}$ such that $x_i = y_{i+k}$ for all $i \in \mathbb
  N$.
\begin{proof} Note that the shift is not injective since we assume
  that $\sigma(S) = S$ and that $S$ is infinite. Combine Theorem
  \ref{matsumoto}, Lemma \ref{sss} and Theorem \ref{ultsimplthm}.
\end{proof}
\end{thm}

Similarly, for subshifts Theorem \ref{relsimplealg} can now be re-formulated as follows:

\begin{thm}\label{subshiftRsimple} Let $S$ be an infinite
  one-sided subshift. Then $C^*_r\left(R_S\right)$ is simple if and
  only if
  the following holds: For every non-empty generalized cylinder
  $C(u;F)$ there is an $m \in \mathbb N$ with the property that for
  all $x \in S$ there is an element $y \in C(u;F)$ such that $x_i = y_{i+m}$ for all $i \in \mathbb
  N$.
\end{thm}

Concerning the existence of a Cartan subalgebra of $\mathcal O_S$ note
that a subshift only has finitely many periodic points of each
period. We can therefore combine Theorem \ref{matsumoto} and Theorem \ref{cartan13} to obtain the following

\begin{thm}\label{cartan2} The Carlsen-Matsumoto algebra $\mathcal O_S$ of
  a subshift $S$ contains a Cartan subalgebra in the sense of Renault,
  \cite{Re3}.
\end{thm}

\end{document}